\documentclass[11pt,twoside]{amsart}
\usepackage{amscd,redpairs}

\author{Valery Alexeev and Michel Brion}
\address{Department of Mathematics\\
University of Georgia\\
Athens, GA 30602, USA}
\email{valery@math.uga.edu}
\address{Institut Fourier, B. P. 74\\
38402 Saint-Martin d'H\`eres Cedex, France}
\email{Michel.Brion@ujf-grenoble.fr}

\begin{document}

\bibliographystyle{amsalpha}

\title{Stable reductive varieties I: Affine varieties}

\date{September 2, 2003}

\maketitle 

\tableofcontents

\setcounter{section}{-1}

\section{Introduction}

The main motivation for this work is to extend the construction of the
moduli spaces of polarized abelian varieties $A_g$ and their
compactifications to a non-commutative setting. Recall that $A_g$
admits many natural compactifications: toroidal (due to Mumford et
al), Satake-Baily-Borel, Borel-Serre etc. Of these, one special
toroidal compactification has an interpretation as the moduli space of
projective stable pairs $(X,D)$: see \cite{Alexeev_CMAV}, which
extends the much earlier work of Namikawa \cite{Namikawa_NewCompBoth}. 
The former paper also constructs the moduli space of stable pairs with
torus action, and semiabelian action. 

This raises the question of defining a nice class of varieties with an
action of a given non-commutative algebraic group $G$, which includes
the group itself, its equivariant compactifications, and also all of
their natural degenerations, ``stable'' varieties. In addition, one
would like to construct a nice, complete moduli space of these stable
varieties, so to say, a toroidal compactification.

That is exactly what we do here, for a connected reductive group $G$. 
It comes equipped with two commuting $G$-actions by left and right
multiplication, that is, with an action of $\gxg$; this reduces to a
$G$-action if $G$ is commutative. We define a class of
\emph{reductive varieties}, consisting of (irreducible) varieties with a
$\gxg$-action. It contains the group $G$ and all of its equivariant
embeddings, and is closed under (flat) degenerations with reduced
and irreducible fibers. We also introduce a larger class of (possibly
reducible) \emph{stable reductive varieties}, closed under 
degenerations with reduced fibers. 

In the commutative case, where $G$ is a torus, reductive varieties are
just toric varieties; in general, they are certain, highly symmetric,
spherical varieties, and stable reductive varieties are obtained from
them by a simple glueing procedure.

In the present paper, we concentrate on \emph{affine} stable
reductive varieties, the projective case being the topic of our
subsequent paper \cite{Alexeev_BrionII}. We classify them all in terms
of combinatorial data; these turn out to be complexes of cones in the
weight space of $G$, invariant under the action of the Weyl group
$W$. We describe all their one-parameter degenerations. Like in
the toric case, they arise from certain stable reductive varieties for
the larger group $\bG_m\times G$. 

An interesting consequence of the classification is the ``toric
correspondence'': it turns out that stable reductive varieties are in
bijection with stable toric varieties equipped with a compatible
$W$-action.

Another nice feature of reductive varieties is the existence of
an associative multiplication law, making each of them into an
algebraic semigroup (generally without identity element). We
characterize the class of semigroups arising in this way, that we call 
\emph{reductive semigroups}. It contains the class of (normal)
\emph{reductive monoids} studied by Putcha, Renner and Rittatore
(see \cite{Rittatore98} and the survey \cite{Solomon99}), that is, of
algebraic semigroups with an identity element and a reductive unit
group.

Finally, we begin the study of families of stable reductive
varieties and, more generally, of varieties with reductive group
action. Examples of these already occur in work of Popov 
\cite{Popov86} constructing a special degeneration of any affine
$G$-variety, and of Vinberg \cite{Vinberg95} on a remarkable family
of reductive monoids. Given a reductive variety $X$, we construct 
a family \`a la Vinberg of such varieties, with general fiber $X$, and
we prove that any (reduced and irreducible) degeneration of $X$ arises
in that way.

We also construct a fine moduli space for families of $G$-subvarieties
of a fixed $G$-module. The construction applies to all
multiplicity-finite $G$-varieties (not only to stable reductive
varieties), and is an application and a generalization of
the multigraded Hilbert scheme of Haiman and Sturmfels 
\cite{HaimanSturmfels}. This result plays a fundamental r\^ole in our
follow-up paper \cite{Alexeev_BrionII}; there we construct the moduli
space of projective pairs $(X,D)$ analogous to the
Deligne-Mumford-Knudsen moduli space of stable curves, and explore
connections with log Minimal Model Program.


\section{Main definitions and results}

We will use \cite{Hartshorne77} as a general reference for algebraic
geometry, and \cite{Popov_Vinberg94}, \cite{Grosshans97} for algebraic
transformation groups.

\begin{notation} (on groups).

The ground field $k$ is algebraically closed, of characteristic zero.
Let $G$ be a connected reductive group over $k$. Let $B$, $B^-$ be
opposite Borel subgroups of $G$, with unipotent radicals $U$, $U^-$.
Then $T=B\cap B^-$ is a maximal torus of $G$; let $N=N_G(T)$ be its
normalizer, and $W=N_G(T)/T$ the corresponding Weyl group. There
exists a unique automorphism $\theta$ of $G$ such that
$\theta(t)=t^{-1}$ for any $t\in T$. Then $\theta$ is involutive; it
exchanges $B$ with $B^-$, stabilizes $N$ and $T$, and fixes pointwise
$W$.

The product $\gxg$ is a connected reductive group; the subgroups
$B^-\times B$, $B\times B^-$ are opposite Borel subgroups, with common
torus $T\times T$. The map 
$\Theta:(g_1,g_2)\mapsto (\theta(g_2),\theta(g_1))$ is an involutive
automorphism of $\gxg$ (note the additional twist). The diagonal
subgroups $\diag G$, $\diag N$, $\diag T$ are
$\Theta$-invariant, and the induced map on $\diag W$ is the identity.
\end{notation}

\begin{notation} (on varieties).
By a \emph{variety}, we mean a reduced scheme of finite type over $k$;
in particular, varieties need not be irreducible. A \emph{$G$-variety}
$X$ is a variety endowed with an action of $G$. The group of
equivariant automorphisms of $X$ is then denoted by $\Aut^G(X)$. As an
example, for any homogeneous space $G/H$, where $H$ is a closed
subgroup of $G$, we have $\Aut^G(G/H) = N_G(H)/H$. The 
fixed point subset of $H$ in $X$, regarded as a closed reduced
subscheme, is denoted by $X^H$; it has an action of $N_G(H)/H$. The
stabilizer of a closed subset $Y\subseteq X$ is the closed subgroup
$\Stab_G(Y)=\{g\in G~\vert~ gY=Y\}$ of $G$.

We will mostly consider $\gxg$-varieties. The \emph{adjoint} of such
a variety $X$ is the $\gxg$-variety $\Theta(X)$ with underlying
variety $X$, and action twisted by $\Theta$. (If $X$ is a $G$-variety,
we define another $G$-variety $\theta(X)$ in a similar way.)
A $G\times G$-variety $X$ is \emph{self-adjoint} if it admits an
involutive automorphism $\Theta_X$ such that 
$\Theta_X(\gamma x)=\Theta(\gamma)\Theta_X(x)$ 
for all $\gamma\in\gxg$ and $x\in X$. As an easy example, any
$T$-variety $X$ is a self-adjoint $T\times T$-variety for the
action $(t_1,t_2)\cdot x=t_1 t_2^{-1}x$, and the identity
automorphism.

For any $\gxg$-variety $X$, the fixed point subset $X^{\diag T}$ is
equipped with an action of 
$N_{\gxg}(\diag T)/\diag T = (\diag N)(T\times T)/\diag T$. Via
the second projection, this group is isomorphic to $WT$, 
the semi-direct product of $T$ with $W\subset \Aut T$.

Finally, recall that a variety $Y$ is \emph{seminormal} if every
finite bijective morphism from a variety to $Y$ is an isomorphism. 
Any variety $X$ admits a unique (hence, equivariant) seminormalization 
$\pi:Y\to X$.

\end{notation}

\begin{definition}\label{stablereductive}
An \emph{affine stable reductive variety} (resp.~an 
\emph{affine reductive variety}) for $G$ is a connected
(resp.~irreducible) affine $\gxg$-variety $X$ satisfying the
following conditions: 

\begin{enumerate}

\item (on singularities) $X$ is seminormal,

\item (on stabilizers) for any $x\in X$, the stabilizer
$\Stab_{\gxg}(x)$ is connected,

\item (on orbits) $X$ contains only finitely many $\gxg$-orbits, 

\item (group-like condition) $(U^-\times U)X^{\diag T}$ contains a
dense subset of every $\gxg$-orbit.

\end{enumerate}

\end{definition}
Here we are following a well-established (but somewhat confusing)
tradition started by Deligne and Mumford \cite{DeligneMumford69} with
their definition of stable curves. Namely, an affine reductive variety
is also a stable reductive variety but not necessarily vice versa. In
other words, ``stable'' serves to widen the class, not to narrow it down.
And since most varieties considered in this paper are affine, we will
often drop the adjective ``affine'', and deal with (stable) reductive
varieties. 

The easiest reductive variety is the group $G$ itself, 
where $\gxg$ acts by $(g_1,g_2)\cdot g = g_1 g g_2^{-1}$: then
$G\simeq (\gxg)/\diag G$. Moreover, $G^{\diag T} = T$, and the subset
$(U^-\times U)G^{\diag T}= U^-TU=B^-B$ 
is open in $G$. Note that $G$ is self-adjoint for
$\Theta_G:g\mapsto\theta(g^{-1})$; moreover, $\Theta_G$ fixes $T$
pointwise.

If $G=T$ is a torus, then condition (4) just means that $\diag T$
fixes pointwise $X$. Therefore, the stable reductive varieties for
$T$ are exactly the seminormal affine varieties where $T$ acts
with finitely many orbits and connected isotropy groups. These are
the affine \emph{stable toric varieties} in the sense of
\cite{Alexeev_CMAV}.

Our first main result (Theorems \ref{correspondence1} and
\ref{correspondence2}) classifies all stable reductive varieties for
$G$, in terms of stable toric varieties for $T$. 

\medskip

\noindent
{\bf Theorem.}
\emph{Let $X$ be a stable reductive variety for $G$. Then 
$X^{\diag T}$ is a stable toric variety for $T$, with a 
compatible action of $W$ (that is, an action of $WT$). 
Moreover, $X$ is self-adjoint for a unique automorphism $\Theta_X$
fixing pointwise $X^{\diag T}$.\\
Any $\gxg$-orbit $\cO$ in $X$ meets $X^{\diag T}$ into a unique
$WT$-orbit; moreover, $\Aut^{\gxg}(X)\simeq \Aut^{WT}(X^{\diag T})$.\\
The assignment $X\mapsto X^{\diag T}$ defines a bijective
correspondence from the stable reductive varieties (for $G$) to the
stable toric varieties (for $T$) with a compatible $W$-action. 
Moreover, $X$ is irreducible if and only if $X^{\diag T}/W$ is.}

\begin{definition}\label{stablesemigroup}
A \emph{stable reductive semigroup} (resp.~a \emph{reductive semigroup})
for $G$ is a connected (resp.~irreducible) affine $\gxg$-variety $X$ 
satisfying conditions (1), (2) of Definition \ref{stablereductive},
and equipped with a morphism $m:X\times X\to X$ satisfying the
following conditions:

\begin{enumerate}

\item (on semigroups) $m$ is an associative multiplication law on $X$. 
 
\item (on equivariance) 
$m((g_1,g)x_1,(g,g_2)x_2)=(g_1,g_2)m(x_1,x_2)$
for all $g_1,g,g_2\in G$ and $x_1,x_2\in X$.

\item (on invariance) $m$ is a categorical quotient for the
$G$-action given by $g\cdot(x_1,x_2)=((1,g)x_1,(g,1)x_2)$.

\end{enumerate}

If, in addition, $X$ is self-adjoint for an automorphism $\Theta_X$
such that 
$$
\Theta_X(m(x_1,x_2))=m(\Theta_X(x_2),\Theta_X(x_1))
$$ 
for all $x_1,x_2\in X$, then $X$ is a 
\emph{self-adjoint reductive semigroup}. 
\end{definition}

The multiplication of a stable reductive semigroup $X$ will be denoted
by $(x_1,x_2)\mapsto x_1 x_2$, and the $\gxg$-action by
$(g_1,g_2,x)\mapsto g_1 x g_2^{-1}$; this notation makes sense by the
equivariance condition.

By Lemma \ref{monoid}, examples of self-adjoint reductive
semigroups include normal \emph{reductive monoids}, that is, 
irreducible algebraic monoids with unit group $G$, where $\gxg$ acts
by left and right multiplication. These are exactly the normal affine
embeddings of the homogeneous space $(\gxg)/\diag G$ (see
\cite{Rittatore98} Theorem 1 and Proposition 1).

The self-adjoint reductive semigroups turn out to be closely related
to reductive varieties, as shown by our second main result (Lemma
\ref{square} and Theorem \ref{semigroup}).

\medskip

\noindent
{\bf Theorem.}
\emph{Any (stable) self-adjoint reductive semigroup is a (stable)
reductive variety.\\ 
Conversely, any reductive variety admits a structure of self-adjoint
reductive semigroup, and any two such structures are conjugate by a
unique automorphism.}

\medskip

We also characterize those stable reductive varieties that admit a
structure of self-adjoint stable reductive semigroup (Proposition
\ref{stablesemired}). The class of (not necessarily self-adjoint)
stable reductive semigroups is much wider: consider, for example, the
$\gxg$-variety $X=Y\times Z$, where $Y$ and $Z$ are affine
$G$-varieties. Then the map 
$m : X\times X\to X,~(y_1,z_1,y_2,z_2)\mapsto (y_1,z_2)$ satisfies
conditions (1) and (2) of Definition \ref{stablesemigroup}. If, in
addition, the categorical quotient of $X$ by the action of $\diag G$
is a point, then $m$ satisfies (3) as well.

\begin{notation} (on schemes).
By a \emph{scheme}, we mean a separated Noetherian scheme over $k$; 
morphisms (resp.~products) of schemes are understood to be
$k$-morphisms (resp.~products over $k$). For any point $s$ of a scheme
$S$, we denote by $k(s)$ the residue field at $s$. The choice of an
algebraic closure $k(\os)\supseteq k(s)$ defines a geometric point
$\os$ over $s$.
\end{notation}

\begin{definition}\label{family}
A \emph{family of affine (stable) reductive varieties} over a scheme
$S$ is a scheme $\cX$ equipped with a morphism $\pi:\cX\to S$ and with
an action of the constant group scheme $\gxg\times S$ over $S$, 
satisfying the following conditions:

\begin{enumerate}

\item $\pi$ is flat and affine. 

\item The geometric fiber $\cX_{\os}$ at every geometric point
$\os$ is a (stable) reductive variety for $(\gxg)(k(\os))$.

\end{enumerate}

\end{definition}

In that case, $\pi$ is the categorical quotient by $\gxg$ (Lemma
\ref{decomposition}). There are obvious notions of trivial
(resp.~locally trivial) families, and one sees easily that any family
of toric varieties is locally trivial (Lemma \ref{loctriv}). Clearly,
this does not extend to families of reductive varieties. But any such
family $\pi:\cX\to S$, where $S$ is integral, becomes locally trivial
after some \'etale base change (Corollary \ref{isotrivial}); this
defines the \emph{general fiber} of $\pi$.

Our third main result concerns moduli. For a fixed $G$-module $V$, we
construct a fine moduli space $\cM_{h,V}$ of $G$-subvarieties of $V$
that have a fixed Hilbert function $h$, and we prove that $\cM_{h,V}$
is quasiprojective (Theorem~\ref{Hilb}). We also describe all
one-parameter flat degenerations of stable reductive varieties
(Proposition~\ref{prop:curve-stable}).

Our fourth main result (Theorem \ref{locstruc}) provides a local model
for families of reductive varieties over an integral scheme, with a
prescribed general fiber. Its starting point is the construction
of a remarkable class of such families, by Vinberg \cite{Vinberg95}. 
To each semisimple group $G_0$, he associated its ``enveloping monoid''
$\Env(G_0)$. This is a normal reductive monoid having $G_0$ as the
derived subgroup of its unit group, and such that the categorical
quotient $\pi:\Env(G_0)\to \Env(G_0)//(G_0\times G_0)$ is flat, with
reduced and irreducible geometric fibers. Moreover, any reductive
monoid satisfying these properties is obtained from $\Env(G_0)$ by 
base change. In addition, $\Env(G_0)//(G_0\times G_0)$ is the affine
space $\bA^r$, where $r$ is the rank of $G_0$; and the fiber of $\pi$
at the identity is $G_0$. Hence $\pi$ is a family of reductive
varieties for $G_0$, with general fiber $G_0$ itself.

Given any reductive variety $X$ for $G$, a variant of Vinberg's
construction yields a family of reductive varieties
$\pi_X:\cV_X\to\bA^r$ with general fiber $X$, the \emph{Vinberg family} 
of $X$; here $r$ is the semisimple rank of $G$ (see 7.5 for details).

\medskip

\noindent
{\bf Theorem.}
\emph{Any family of reductive varieties over an integral scheme is
locally a base change of the Vinberg family of its general fiber.}


\section{General criteria}

\label{sec: General criteria}


\subsection{Seminormality and connectedness of isotropy groups}

\label{subsection: Seminormality and connectedness of isotropy groups}

We will obtain criteria for conditions (1), (2) of Definition
\ref{stablereductive} to hold; we begin with some additional notation
on groups and varieties. 

The weight lattice of $G$ is the character group of $T$, denoted by
$\Lambda$; let $\Lambda_{\bR}$ be the corresponding real vector space.
The positive Weyl chamber of $\Lambda_{\bR}$ associated with $B$ is
denoted by $\Lambda^+_{\bR}$; then $\Lambda^+=\Lambda\cap
\Lambda^+_{\bR}$ is the set of dominant weights. For each
$\lambda\in\Lambda^+$, let $V_\lambda$ be the simple $G$-module with
highest weight $\lambda$.

For any $G$-variety $X$, the algebra of regular functions
$k[X]$ is a rational $G$-module; thus, it decomposes as a direct
sum of simple submodules $V_{\lambda}$, with (possibly infinite)
positive multiplicities. The set of all such $\lambda$ is the
\emph{weight set} $\Lambda^+_X$; let $\Lambda_X$ (resp.~$C_X$)
be the subgroup of $\Lambda$ (resp.~the cone of $\Lambda_{\bR}$)
generated by $\Lambda^+_X$. Then 
$\Lambda^+_X\subseteq \Lambda^+\cap C_X$; this inclusion may be
strict.

If, in addition, $X$ is affine, then the invariant subalgebra
$k[X]^G$ is finitely generated. We denote by $X//G$ the
corresponding affine variety; the corresponding morphism  
$p:X\to X//G$ is the categorical quotient by $G$. The invariant
subalgebra $k[X]^U$ is finitely generated as well; moreover, $T$
acts on $k[X]^U$ and hence on the associated variety $X//U$. For any
$\lambda\in\Lambda$, let $k[X]^U_\lambda$ be the $\lambda$-weight
space in the $T$-module $k[X]^U$. Then 
$k[X]^U_\lambda \simeq \Hom^G(V_\lambda,k[X])$; in particular, 
$k[X]^U_0=k[X]^G$. Any $k[X]^U_\lambda$ is a finitely generated
$k[X]^G$-module; its dimension (as a $k$-vector space) is the
multiplicity of the $G$-module $V_\lambda$ in $k[X]$ if
$\lambda\in\Lambda^+$, and is $0$ otherwise. 

We say that the affine $G$-variety $X$ is \emph{multiplicity-free}
(resp.~\emph{multiplicity-finite, multiplicity-bounded}) if all
these multiplicities are $0$ or $1$ (resp.~finite, bounded), see
\cite{Grosshans97} \S 11. Note that $X$ is multiplicity-finite if and
only if $k[X]^G$ is finite-dimensional as a $k$-vector space.

Recall that a $G$-variety $X$ is \emph{spherical} if $X$ is normal
and contains a dense $B$-orbit; then it contains only finitely many
$B$-orbits, see \cite{Grosshans97} Theorem 22.6. The affine spherical
varieties are exactly the irreducible normal, multiplicity-free
varieties, by \cite{Grosshans97} Theorem 11.1.

\begin{lemma}\label{Uinv}

Let $X$ be an affine $G$-variety. 

\begin{enumerate}

\item $X$ is connected if and only if $X//G$ is. As a consequence, any
  multiplicity-free variety is connected.

\item If $X$ is irreducible (resp.~normal, seminormal), then so is $X//G$. 

\item $X$ is irreducible (resp.~normal, seminormal) if and only if $X//U$ is.

\end{enumerate}

\end{lemma}

\begin{proof}
The assertions on connectedness, irreducibility, and normality are
well-known, see e.g. \cite{Grosshans97} \S 18. 

(2) Assume that $X$ is seminormal. Let $\pi:Y\to X//G$ be the
    seminormalization. Let $p:X\to X//G$ be the quotient morphism, 
    and form the cartesian square
$$
\CD
Z @>{q}>>  Y \\
@VVV @V{\pi}VV \\
X @>{p}>> X//G. \\
\endCD
$$
Let $Z_{\red}$ be the reduced subscheme of $Z$. Then the induced map
$Z_{\red}\to X$ is finite and bijective, hence an
isomorphism. Moreover, $G$ acts on $Z$, and $q:Z\to Y$ is the quotient
map (since $G$ is reductive). Thus, this map factors through an
isomorphism $Z_{\red}//G\to Y$ (since $Y$ is reduced). It follows that
$\pi:Y\to X//G$ is an isomorphism as well.

(3) If $X$ is seminormal, then so is $X\times G//U$ (since
    $G//U$ is normal.) Thus, $(X\times G//U)//G$ is seminormal
    by (1). But $(X\times G//U)//G\simeq X//U$ by \cite{Grosshans97} 
    Theorem 9.1, so that $X//U$ is seminormal.

Conversely, assume that $X//U$ is seminormal and consider the
seminormalization $\pi:Y\to X$. Then $G$ acts on $Y$ and the map 
$\pi\times \id:Y\times G//U\to X\times G//U$ is the
seminormalization. It follows that the induced map 
$(Y\times G//U)//G\to (X\times G//U)//G$ is finite and bijective. 
By \cite{Grosshans97} Theorem 9.1 again, this means that the map 
$\pi//U: Y//U\to X//U$ is finite and bijective. Since $X//U$ is
seminormal, $\pi//U$ is an isomorphism; since $Y$ is a $G$-variety,
it follows that $\pi$ is an isomorphism. 
\end{proof}

\begin{lemma}\label{multfree}
Let $X$ be an affine spherical $G$-variety. 
Then: 

\begin{enumerate}

\item $C_X$ is a rational polyhedral convex cone, and
  $\Lambda_X^+=C_X\cap\Lambda_X$. 

\item Any closed irreducible $G$-subvariety $Y$ of $X$ is uniquely
  determined by $\Lambda_Y^+$. Moreover, $Y$ is normal, and there
  exists a unique face $F$ of $C_X$ such that
  $\Lambda_Y^+=\Lambda_X^+\cap F$. 

\item If $\Lambda_X^+$ is saturated in $\Lambda$, then each isotropy
  group of $G$ in $X$ is connected.

\end{enumerate}

\end{lemma}

\begin{proof}
Note that $X//U$ is an affine toric variety for the torus with
character group $\Lambda_X$ (a quotient of $T$). By the theory of
toric varieties (see e.g. \cite{Fulton93}), this implies (1), and
(by considering $Y//U$) the second and third assertions of (2). Since
$X$ is multiplicity-free and the ideal of $Y$ in $k[X]$ is
$G$-invariant, this ideal is uniquely determined by
$\Lambda_Y^+$. This completes the proof of (2).

(3) By (2), $\Lambda_Y^+$ is saturated in $\Lambda$ for any closed
    irreducible $G$-subvariety $Y$. Thus, it suffices to show that
    the isotropy group $H$ of a point of the open $G$-orbit is
    connected. Let $H^0$ be the connected component of the identity in
    $H$. Then we have a dominant $G$-equivariant map $G/H^0\to X$. 
    Let $\tilde X$ be the normalization of $X$ in 
    the field of rational functions on $G/H^0$. This is an affine
    spherical $G$-variety, and the algebra $k[\tilde X]^U$ is
    integral over $k[X]^U$. It follows that any
    $\chi\in\Lambda_{\tilde X}^+$ admits a positive integral
    multiple in $\Lambda_X^+$. Since $\Lambda_X^+$ is saturated in
    $\Lambda$, this implies $\Lambda_{\tilde X}^+=\Lambda_X^+$, whence
    $\tilde X=X$, and $H^0=H$.
\end{proof}

\begin{lemma}\label{saturated}
Let $X$ be an affine $T$-variety. Then $X$ is multiplicity-bounded if
and only if every irreducible component is multiplicity-free; in that
case, $X$ contains only finitely many $T$-orbits. 

If, in addition, $X$ is multiplicity-free, then the following
conditions are equivalent:

\begin{enumerate}

\item $X$ is seminormal, with connected isotropy groups.

\item $\Lambda_X^+$ is saturated in $\Lambda$.

\end{enumerate}

Then any $T$-orbit closure in $X$ is normal. 
\end{lemma}

\begin{proof}
If every irreducible component of $X$ is multiplicity-free, then $X$
is clearly multiplicity-bounded. Conversely, let $X$ be
multiplicity-bounded and let $Y$ be an irreducible component. If we
can find linearly independent $T$-eigenvectors $\varphi,\psi\in k[Y]$
of the same weight $\lambda$, then for any positive integer $n$, the
monomials 
$\varphi^n,\varphi^{n-1}\psi,\ldots,\varphi\psi^{n-1},\psi^n$ 
are linearly independent $T$-eigenvectors of weight $n\lambda$. This
contradicts the assumption that $X$ is multiplicity-bounded; thus,
$Y$ is multiplicity-free. Using Lemma \ref{multfree}, it follows that
$X$ contains only finitely many orbits.

(1) $\Rightarrow$ (2): Let $(X_i)$ be the (finite) family of all closed 
    irreducible $T$-subvarieties of $X$, and let $(\tilde X_i)$ be
    the family of their normalizations. These form a direct system,
    with direct limit $X$ by \cite{Alexeev_CMAV} Theorem 2.3.14. In
    other words, $k[X]=\varprojlim_i k[\tilde X_i]$. Since each 
    $\tilde X_i$ is a toric variety for a quotient of $T$ by a
    subtorus, $\Lambda_{\tilde X_i}^+$ is saturated in $\Lambda$. It
    follows that the same holds for $\Lambda_X^+$.

(2) $\Rightarrow$ (1): Let $Y$ be a closed irreducible $T$-subvariety 
    of $X$. Then the assumptions imply that $Y$ is normal 
    and multiplicity-free, and that $\Lambda_Y^+$ is saturated in
    $\Lambda$. Thus, $Y$ is a toric variety for the quotient of
    $T$ by a subtorus. So all isotropy subgroups of $X$ are connected. 
    
    Let $Y$ be the direct limit of all closed irreducible
    $T$-subvarieties of $X$. Then $Y$ is seminormal by
    \cite{Alexeev_CMAV} Corollary 2.3.10. The natural morphism
    $\pi:Y\to X$ is clearly finite and bijective; since $X$ is
    multiplicity-free and $\Lambda_X^+$ is saturated, $\pi$ is an
    isomorphism. In other words, $X$ is seminormal. Moreover,
    every $T$-orbit closure is normal.
\end{proof}

Combining Lemmas \ref{Uinv}, \ref{multfree} and \ref{saturated}, we
obtain the following criterion for an affine $G$-variety to satisfy
conditions (1), (2) and (3) of Definition \ref{stablereductive}.

\begin{proposition}\label{connected}
Any multiplicity-bounded $G$-variety $X$ contains only fini\-tely many
$G$-orbits, and these are spherical. As a consequence, $X$ contains
only finitely many $B$-orbits.

If, in addition, $X$ is multiplicity-free and $\Lambda_X^+$ is
saturated in $\Lambda$, then $X$ is connected and seminormal, the
$G$-orbit closures are normal, and their isotropy groups are
connected. 
\end{proposition}

We will also need the following observation.

\begin{lemma}\label{normality}
Let $X$ be an irreducible, multiplicity-free $G$-variety and let
$\nu:\tilde X \to X$ be the normalization map. If all isotropy groups
are connected, then $\nu$ is bijective. 
\end{lemma}

\begin{proof}
Consider an orbit $\cO$ in $X$, and an orbit $\tilde\cO$ in
$\nu^{-1}(\cO)$. Then the restriction $\nu:\tilde\cO\to \cO$ is an
isomorphism, since the isotropy groups of $\cO$ are connected. Hence
$\nu^{-1}(\cO)$ is a union of finitely many copies of $\cO$.

On the other hand, for any orbit $\tilde\cO$ in $\tilde X$, with
closure $Y$, we have $\Lambda_{\tilde\cO}=\Lambda_Y$. By Lemma
\ref{multfree}, it follows that the linear span of
$\Lambda_{\tilde\cO}$ intersects $C_X$ along a face, which uniquely
determines $Y$ and hence $\tilde\cO$. Thus, any two distinct orbits in
$\tilde X$ are non-isomorphic as $G$-varieties. It follows that
$\nu$ is bijective. 
\end{proof}


\subsection{Finiteness of number of orbits and group-like condition}

\label{subsection:Finiteness of number of orbits and 
group-like condition}

We will obtain a representation-theoretic criterion for a
variety to satisfy conditions (3) and (4) of Definition
\ref{stablereductive}. For this, we introduce some notation on
representations.

The simple $\gxg$-modules are the $V_{\lambda}\otimes_k V_{\mu}$ where
$\lambda,\mu\in\Lambda^+$; those containing non-zero fixed points of
$\diag G$ are the $\End V_\lambda=V_\lambda^*\otimes_k V_\lambda$. Note
that $\theta(V_\lambda)\simeq V_\lambda^*$, so that 
$\Theta(V_{\lambda}\otimes_k V_{\mu})\simeq V_\mu^*\otimes_k V_\lambda^*$;
in particular, $\Theta(\End V_\lambda)\simeq\End V_\lambda$ as
$\gxg$-modules.

\begin{definition}\label{diagonal}
An affine $\gxg$-variety $X$ is \emph{diagonal} if every simple
submodule of the $\gxg$-module $k[X]$ is isomorphic to some 
$\End V_\lambda$. Equivalently, the weights of the $T\times T$-module 
$k[X]^{U^-\times U}$ are of the form $(-\lambda,\lambda)$, where
$\lambda\in\Lambda^+$. We then identify the weight set of $X$ to a
subset of $\Lambda^+$, via the second projection; this identifies
$\Lambda_X$ to a subgroup of $\Lambda$.
\end{definition}

\begin{proposition}\label{module}
For an affine $\gxg$-variety $X$, the following conditions 
are equivalent:

\begin{enumerate}

\item $X$ contains only finitely many $\gxg$-orbits, and 
$(U^-\times U) X^{\diag T}$ contains a dense subset of any such
orbit.

\item $X$ is diagonal and multiplicity-bounded.

\item $X$ contains only finitely many $B^-\times B$-orbits, and any
such orbit which is open in its $\gxg$-orbit contains fixed points of
$\diag T$.

\end{enumerate}

\end{proposition}

\begin{proof}
(1)$\Rightarrow$(2) Every $T\times T$-eigenvector
$f\in k[X]^{U^-\times U}$ is uniquely determined by its restriction
to $X^{\diag T}$. Thus, the weight of $f$ is $(-\lambda,\lambda)$ for
some $\lambda\in\Lambda^+$. 
Moreover, for any $\gxg$-orbit $\cO$, the subset $\cO^{\diag T}$
contains only finitely many $T\times T$-orbits 
(since $T\times T$ is the centralizer of $\diag T$ in $G\times G$).
As a consequence, $X^{\diag T}$ contains only finitely many 
$T\times T$-orbits. Hence the multiplicities of all 
$T\times T$-weight spaces in $k[X^{\diag T}]$ are bounded, and the
same holds for the multiplicities of all weight spaces
in $k[X]^{U^-\times U}$.

(2)$\Rightarrow$(3) Let $Y$ be a closed irreducible $\gxg$-subvariety 
of $X$. Then $Y$ is multiplicity-free, by Lemma \ref{saturated}. 
It follows that $Y$ contains only finitely
many $B^-\times B$-orbits, and hence an open
$B^-\times B$-orbit $Y_0$. Moreover, the weights of the 
$T\times T$-eigenvectors in $k[Y]^{U^-\times U}$, and hence in 
$k[Y_0]^{U^-\times U}$, are of the form $(-\lambda,\lambda)$. Thus,
$\diag T$ fixes pointwise the quotient $Y_0/(U^-\times U)$. This
implies easily that $Y_0$ contains a $\diag T$-fixed point.

(3)$\Rightarrow$(1) is obvious.
\end{proof} 

The local structure of a variety satisfying these conditions
is described by the following

\begin{lemma}\label{local}
Let $X$ be a diagonal, multiplicity-bounded $\gxg$-variety.

\begin{enumerate}

\item
Let $f\in k[X]^{U^-\times U}$ be a $T\times T$-eigenvector of weight
$\lambda$. Then the stabilizer of the line $k f$ is a product 
$P^-\times P$, where $P\supseteq B$ and $P^-\supseteq B^-$ are
opposite parabolic subgroups. Let $L=P\cap P^-\supseteq T$ be their
common Levi subgroup. Then there exists a locally closed affine
$L\times L$-subvariety $Z\subseteq X$ such that the map
$$
R_u(P^-)\times R_u(P)\times Z \to X,~
(g_1,g_2,x)\mapsto (g_1,g_2)x
$$
is an open immersion with image the principal open subset
$X_f$. Moreover, $Z$ is a diagonal, multiplicity-bounded 
$L\times L$-variety. 

\item
If, in addition, $X$ is multiplicity-free, then so is $Z$.

\item
For any open $B^-\times B$-orbit $X_0$ in $X$, there exists $f$ as in
(1) such that $X_f =X_0$. Then $P$ and $L$ depend only on $X_0$, and 
$Z=X_0^{\diag T}$ is a unique $T\times T$-orbit, fixed pointwise by
$[L,L]\times [L,L]$, where $[L,L]$ is the derived subgroup of $L$.

\end{enumerate}

\end{lemma}

\begin{proof}
(1) The first assertion follows from \cite{Knop94} Theorem 2.3. 
Together with Proposition \ref{module}, it implies the second
assertion.

(2) Let $u$, $v$ be $T\times T$-eigenvectors in 
$k[Z]^{(U^-\cap L)\times(U\cap L)}$, of the same weight $\mu$. 
We may regard them as $T\times T$-eigenvectors in 
$k[X_f]^{U^-\times U}$. Thus, there exists a
positive integer $n$ such that $uf^n$, $vf^n$ lift to 
$T\times T$-eigenvectors $\tilde u$, $\tilde v$ in 
$k[X]^{U^-\times U}$, under the localization
$k[X]\to k[X_f]$. Then $\tilde u$, $\tilde v$ have 
weight $\mu + n \lambda$. Since $X$ is multiplicity-free,
$\tilde v\in k^*\tilde u$, so that $v\in k^* u$.

(3) We may choose $f$ as in (1) such that $f$ vanishes everywhere on
$X-X_0$; then $X_f=X_0$. By \cite{Knop94} Proposition 2.4, the
corresponding $Z$ is fixed pointwise by  $[L,L]\times [L,L]$. Since
$X_0$ is a unique $B^-\times B$-orbit, 
isomorphic to $R_u(P^-)\times R_u(P)\times Z$, it follows that $Z$ is
a unique $T\times T$-orbit, containing $X_0^{\diag T}$. The latter is
$T\times T$-invariant and non-empty; thus, $Z=X_0^{\diag T}$.
\end{proof}


\section{Orbits in stable reductive varieties}

\label{sec: structure}


\subsection{Isotropy groups}

\label{subsec: Isotropy groups}

We will determine those quasi-affine homogeneous 
$G\times G$-varieties $\cO$ that satisfy the ``group-like'' condition
(4) of Definition \ref{stablereductive}. By Proposition 2.7, we may
find $x\in\cO^{\diag T}$ such that $(B^-\times B)x$ is open in
$\cO$. Then the isotropy group $H=\Stab_{\gxg}(x)$ contains $\diag T$,
and the product $(B^-\times B)H$ is open in $\gxg$. To describe $H$,
we need additional notation on subgroups of $G$. 

Let $\Phi\subset\Lambda$ be the root system of $(G,T)$, with its subset
of positive roots $\Phi^+$ (the roots of $(B,T)$), and corresponding
subset of simple roots $\Pi$. The coroot of each $\alpha\in\Phi$ is
denoted by $\check\alpha$. For any subset $I\subseteq \Pi$, let
$P_I\supseteq B$ be the corresponding parabolic subgroup of $G$, and
let $P_I^-\supseteq B^-$ be the opposite parabolic subgroup; then
$L_I=P_I\cap P_I^-$ is a Levi subgroup of both, containing $T$. Let
$W_I\subseteq W$ (resp.~$\Phi_I\subseteq \Phi$) be the Weyl group
(resp.~the root system) of $(L_I,T)$ and let $Z(L_I)=Z_I$
be the center of $L_I$. Then $Z_I\subseteq T$ is the intersection of
the kernels of the characters $\alpha\in I$.

\begin{proposition}\label{structure}
Let $H\subseteq \gxg$ be a closed subgroup satisfying the following
conditions. 

\begin{enumerate}

\item $H$ contains $\diag T$.

\item $(B^-\times B)H$ is open in $\gxg$.

\item The homogeneous space $\gxg/H$ is quasi-affine.

\end{enumerate}
Then there  exist a subset $I\subseteq \Pi$, union of two orthogonal
subsets $J$ and $K$, and a subgroup $H_J\subseteq G$ such that:

\begin{enumerate}

\item  $[L_J,L_J]\subseteq H_J\subseteq [L_J,L_J] Z_K$, and 

\item $H$ is conjugate in $T\times T$ to 
$R_u(P_I) H_J\times R_u(P_I^-)H_J) \diag L_K$.

\end{enumerate}
Moreover, $N_{\gxg}(H)=H(Z_K\times Z_K)$.

\end{proposition}

\begin{proof}
Let $P\supseteq B$ be as in Lemma \ref{local} (3). By that Lemma,
the multiplication map 
$R_u(P^-)\times R_u(P)\times (L\times L)H\to \gxg$
is an open immersion; moreover, $(L\times L)\cap H$ contains
$([L,L]\times [L,L])\diag T$. Write $P=P_J$ for some 
$J\subseteq\Pi$; then $(L_J\times L_J)\cap H=(H_J\times H_J)\diag T$
for a unique closed subgroup $H_J$ of $L_J$, containing $[L_J,L_J]$.
Therefore, one obtains a decomposition of the Lie algebra of $\gxg$:
$$
\Lie(\gxg)=
\Lie(R_u(P_J^-)\times R_u(P_J))\bigoplus(\Lie(L_J\times L_J)+\Lie(H)),
$$
and $\Lie(L_J\times L_J)\cap\Lie(H)=\Lie(H_J\times H_J)+\diag\Lie(T)$. 
On the other hand, one has the decomposition
$$
\Lie(\gxg)=\Lie(R_u(P_J^-)\times R_u(P_J))\oplus\Lie(L_J\times L_J)
\oplus\Lie(R_u(P_J)\times R_u(P_J^-)).
$$

For any $\alpha\in\Phi$, choose a root vector 
$e_{\alpha}\in\Lie(G)_{\alpha}$. Since $\Lie(H)$ is stable
under $\diag(T)$, it decomposes as a sum of weight spaces; thus,
$\Lie(H)$ has a basis consisting of elements of $\Lie(T\times T)$ 
and of certain $(u_{\alpha}e_{\alpha},v_{\alpha}e_{\alpha})$
where $u_{\alpha}$, $v_{\alpha}$ are scalars. In fact, the
preceding decompositions imply that
$$
\Lie(H)=\Lie((L_J\times L_J)\cap H)\oplus
\bigoplus_{\alpha\in\Phi^+ -\Phi_J} 
k(e_{\alpha},c_{\alpha}e_{\alpha})\oplus
\bigoplus_{\alpha\in\Phi^+ -\Phi_J} 
k(d_{\alpha}e_{-\alpha},e_{-\alpha})
$$
for certain scalars $c_{\alpha}$, $d_{\alpha}$. 

We claim that $c_{\alpha}=d_{\alpha}$ for any 
$\alpha\in\Phi^+ - \Phi_I$. To see this, set 
$[e_{\alpha},e_{-\alpha}]=h_{\alpha}$; then 
$h_{\alpha}\in\Lie(T)$, and $[h_{\alpha},e_{\alpha}]$
is a non-zero scalar multiple of $e_{\alpha}$. One has
$$
[(e_{\alpha},c_{\alpha}e_{\alpha}),
(d_{\alpha}e_{-\alpha},e_{-\alpha})]=
(d_{\alpha}h_{\alpha},c_{\alpha}h_{\alpha})\in \Lie(H).
$$
If $c_{\alpha}\ne d_{\alpha}$, then $\Lie(H)$ contains both 
$(h_{\alpha},0)$ and $(0,h_{\alpha})$, since $\Lie(H)$
contains $\diag\Lie(T)$. Thus, $\Lie(H)$ contains
$[(h_{\alpha},0),(e_{\alpha},c_{\alpha}e_{\alpha})]=
([h_{\alpha},e_{\alpha}],0)$, that is, 
$(e_{\alpha},0)\in\Lie(H)$. So $c_{\alpha}=0$; likewise, 
$d_{\alpha}=0$, a contradiction. The claim is proved.

Next we consider the images $H_1$, $H_2$ of $H$ under the 
projections $p_1,p_2:\gxg\to G$. By the decomposition of $\Lie(H)$
into $\diag T$-eigenspaces and the preceding claim, $H_1$ is a
parabolic subgroup containing $B$, and the set of roots of $(H_1,T)$
equals
$$
\Phi^+ \cup \Phi_J \cup 
\{-\alpha~\vert~ \alpha\in \Phi^+ - \Phi_J, ~ c_{\alpha}\ne 0\}.
$$
Moreover, $H_2$ is the opposite parabolic subgroup to $H_1$ that
contains $T$. Hence $H_1=P_I$ and $H_2=P_I^-$ for some subset
$I\subseteq\Pi$ containing $J$; then
$$
\Phi_I^+=\Phi_J^+\cup
\{ \alpha\in \Phi^+ - \Phi_J ~\vert~ c_{\alpha}\ne 0\}.
$$
Let
$E=\{ \alpha\in \Phi^+ - \Phi_J ~\vert~ c_{\alpha}\ne 0\}$. 
If $\alpha\in E$ and $\beta\in\Phi_J^+$, then
$\alpha+\beta\notin\Phi_I^+$ (otherwise, since 
$[(e_{\alpha},c_{\alpha} e_{\alpha}),(e_{\beta},0)]\in\Lie(H)$,
we obtain $(e_{\alpha+\beta},0)\in\Lie(H)$; multiplying by
$(e_{-\beta},0)\in\Lie(H)$, it follows that 
$(e_{\alpha},0)\in\Lie(H)$, a contradiction). Since $\Phi_I^+$ is the
disjoint union of $\Phi_J^+$ and $E$, there exists a subset
$K\subseteq I$ orthogonal to $J$, such that $E=\Phi^+_K$. Moreover,
$I=J\cup K$.

Together with the decomposition of $\Lie(H)$ into 
$\diag T$-eigenspaces, this implies that the connected component
$H^0$ satisfies
$$
H^0= (R_u(P_I)H_J^0\times R_u(P_I^-)H_J^0) H_K,
$$
where $H_K$ is the closed connected subgroup of $\gxg$ with Lie
algebra spanned by $\Lie(\diag T)$ and by the 
$(e_{\alpha},c_{\alpha}e_{\alpha})$, 
$(c_{\alpha}e_{-\alpha},e_{-\alpha})$ ($\alpha\in\Phi^+_K$).
Then $H_K\subseteq L_K\times L_K$, and both projections $H_K\to L_K$
are isomorphisms: $H_K$ is the graph of an automorphism of the group
$L_K$ fixing $T$ pointwise. Such an automorphism is the conjugation 
by an element of $T$; therefore, we may assume that 
$H^0 = (R_u(P_I)H_J^0\times R_u(P_I^-)H_J^0)\diag L_K$.
Moreover, $H_J^0=[L_J,L_J](T\cap H_J^0)$, and by considering the
adjoint action of $T\cap H_J^0$ on $(e_{\alpha},e_{\alpha})$
where $\alpha\in K$, one obtains that $\alpha$ vanishes on 
$T\cap H_J^0$. Thus, $[L_J,L_J]\subseteq H_J^0\subseteq [L_J,L_J] Z_K$.

Let $\tilde H$ be the normalizer of $H^0$ in $\gxg$. Then $\tilde H$
normalizes the unipotent radical $R_u(P_I)\times R_u(P_I^-)$ 
of $H^0$, so that $\tilde H\subseteq P_I\times P_I^-$.
In fact, 
$\tilde H = (R_u(P_I)\times R_u(P_I^-))(\tilde H\cap(L_I\times L_I))$,
and $\tilde H\cap (L_I\times L_I)$ is the normalizer in 
$L_I\times L_I$ of $(H_J^0\times H_J^0)\diag L_K$. Now the preceding
arguments show that 
$\tilde H = (R_u(P_I) [L_I,L_I] Z_K \times R_u(P_I^-) [L_I,L_I] Z_K) 
\diag L_K$.
Since 
$H^0\subseteq H\subseteq N_{\gxg}(H)\subseteq N_{\gxg}(H^0)=\tilde H$, 
this implies our assertions on the structure of $H$ and its normalizer.
\end{proof}

For $H= (R_u(P_I) H_J\times R_u(P_I^-)H_J) \diag L_K$ as in
Proposition \ref{structure}, let 
$$
T_H = \{t\in T~\vert~ (t,1)\in H\} = T \cap H_J
$$
and denote by $\Lambda_H$ the set of those characters of $T$ that
vanish on $T_H$. Then $K\subset \Lambda_H \subseteq J^{\perp}$, since
the coroots of elements of $J$ are one-parameter subgroups of $T_H$.
We call the triple $(K,\Lambda_H,J)$ the 
\emph{combinatorial invariant} of $H$. It characterizes the orbit
$G\times G/H$, by Proposition \ref{structure}; a criterion for that
orbit to be quasi-affine will be given in Proposition
\ref{quasiaffine} below.

Conversely, given two orthogonal subsets $J,K\subseteq\Pi$ and a
subgroup $\Lambda'\subseteq \Lambda$ such that
$K\subset \Lambda' \subseteq J^{\perp}$, let $T'$ be the subgroup of
$T$, intersection of the kernels of all the $\chi\in\Lambda'$.
Also, let $H_J = [L_J,L_J] T'$ and $I=J\cup K$. Then
$$
H_{K,\Lambda',J}=(R_u(P_I) H_J \times R_u(P_I^-) H_J) \diag L_K
$$
is a subgroup of $G\times G$ containing $\diag T$, such that
$(B^-\times B) H_{K,\Lambda',J}$ is open in $\gxg$, and with
combinatorial invariant $(K,\Lambda',J)$. We record some easy
properties of the corresponding homogeneous spaces. 

\begin{lemma}\label{automorphisms}
For $H=H_{K,\Lambda',J}$, one has:

\begin{enumerate}

\item $H$ is connected $\Leftrightarrow$ $T'$ is connected
$\Leftrightarrow$ $\Lambda'$ is saturated in $\Lambda$.
 
\item $\Aut^{\gxg}(\gxg/H) = N_{\gxg}(H)/H$ is isomorphic to $Z_K/T_H$, a
diagonalizable group with character group $\Lambda'/\bZ K$.

\item $\dim \gxg/H = \rank \Lambda' + \vert \Phi-\Phi_J\vert$.

\item $\Theta(H)=H$, so that $\Theta$ induces an involutive
automorphism $\Theta_{\gxg/H}$ fixing $(T\times T)H/H$ pointwise. 
This makes $\gxg/H$ a self-adjoint $\gxg$-variety. 
\end{enumerate}
\end{lemma}

\begin{proof}
(1) Since $H^0=(R_u(P_I) H_J^0 \times R_u(P_I^-) H_J^0) \diag L_K$,
and $[L_J,L_J]\subseteq H_J^0\subseteq L_J$, one has
$H/H^0 \simeq H_J/H_J^0 \simeq (T\cap H_J)/(T\cap H_J)^0 =T'/T'{}^0$. 
Moreover, $T'$ (resp.~$T'{}^0$) is the intersection of the kernels
of all characters in $\Lambda'$ (resp.~in the saturation of $\Lambda'$
in $\Lambda$).

(2), (3) and (4) are straightforward.
\end{proof}


\subsection{Algebras of regular functions}

\label{subsec: Algebras of regular functions}

We study the algebra $k[\gxg/H]$, where $H=H_{K,\Lambda',J}$. We
begin with the simplest case where $H=\diag G$, that is, $J=\emptyset$,
$\Lambda'=\Lambda$, and $K=\Pi$. As is well-known, the algebra
$k[\gxg/\diag G]\simeq k[\diag G]$ is spanned by matrix
coefficients of simple $G$-modules. Specifically, the direct sum of
the maps
$$
\End V_\lambda\to k[G],~u\mapsto (g\mapsto \Tr_{V_\lambda}(ug))
$$
yields an isomorphism of $\gxg$-modules:
$\bigoplus_{\lambda\in\Lambda^+}\End V_{\lambda}\simeq k[G].$

To describe the multiplication in $k[G]$, consider the decomposition
of tensor products of $G$-modules: 
$$
V_{\lambda}\otimes_k V_{\mu} \simeq 
\bigoplus_{\nu} N^{\nu}_{\lambda\mu} \otimes_k V_{\nu},
$$
where 
$N_{\lambda\mu}^{\nu}=\Hom^G(V_{\nu},V_{\lambda}\otimes_k V_{\mu})$. 
This defines $\gxg$-equivariant projections
$$
p^{\nu}_{\lambda\mu}: \End V_{\lambda} \otimes_k \End V_{\mu} \to 
\End N^{\nu}_{\lambda\mu} \otimes_k \End V_{\nu}.
$$
Then the multiplication in $k[G]$ of any 
$f_{\lambda}\in \End V_{\lambda}$ and $f_{\mu}\in \End V_{\mu}$
is given by
$$
f_{\lambda}\cdot f_{\mu} = \sum_{\nu} 
\Tr_1( p^{\nu}_{\lambda\mu}(f_{\lambda}\otimes f_{\mu})),
$$
where 
$\Tr_1:\End N^{\nu}_{\lambda\mu}\otimes_k\End V_{\nu}\to\End V_{\nu}$ 
is the trace in the first summand.

For any $\lambda\in\Lambda^+$, let 
$\chi_\lambda:g\mapsto \Tr_{V_\lambda}(g)$ be the corresponding
character. Then the $\chi_\lambda$, $\lambda\in \Lambda^+$, are a
basis of the space $k[G]^{\diag G}$ of conjugation-invariant
functions. The multiplication of characters is given by
$$
\chi_\lambda \cdot \chi_\mu = \sum_{\nu} 
c_{\lambda\mu}^{\nu} \chi_\nu,
$$
where the $c_{\lambda\mu}^{\nu}=\dim N_{\lambda\mu}^{\nu}$
are the ``Littlewood-Richardson coefficients''. Note that
$c_{\lambda\mu}^{\nu}=0$ unless $\lambda+\mu-\nu$ is a linear
combination of elements of $\Pi$ with nonnegative integer coefficients; 
we write $\nu\leq_{\Pi}\lambda+\mu$. In particular, the product
$\End V_\lambda \cdot \End V_\mu$ is the sum of those $\End V_\nu$ such
that $c_{\lambda\mu}^{\nu}\neq 0$.

We now generalize this to any homogeneous space
$\gxg/H_{K,\Lambda',J}$. For this, we construct analogues of
matrix coefficients.

For any simple $G$-module $V_{\lambda}$, let $V_{\lambda,K}$
be the sum of the $T$-weight spaces in $V_{\lambda}$ associated with
those weights $\chi$ such that $\chi\leq_K\lambda$, and let
$V_{\lambda}^K$ be the sum of all other weight spaces. 
Then $V_{\lambda,K}=V_{\lambda}^{R_u(P_K)}$ is a simple $L_K$-module
with highest weight $\lambda$. Moreover, $V_{\lambda}^K$ is stable
under $P_K^-$, and $V_{\lambda}=V_{\lambda,K}\oplus V_{\lambda}^K$.
Let 
$$
p_{\lambda,K}:V_{\lambda}\to V_{\lambda}
$$
be the projection to $V_{\lambda,K}$ with kernel $V_{\lambda}^K$.

\begin{proposition}\label{product}

Let $H=H_{K,\Lambda',J}$ and $\gxg/H=\cO$.

\begin{enumerate}

\item For any $\lambda\in \Lambda'\cap \Lambda^+$, the projection
$p_{\lambda,K}\in \End V_{\lambda}$ is fixed by $H$. Thus, the map
$$
\End V_{\lambda}\to k[\gxg],~u\mapsto ((x,y)\mapsto
\Tr_{V_{\lambda}}(u x p_{\lambda,K} y^{-1}))
$$
is a $\gxg$-equivariant morphism, with image in $k[\cO]$.

\item The resulting map 
$$
\bigoplus_{\lambda\in\Lambda'\cap \Lambda^+}\End V_{\lambda}
\to k[\cO]
$$
is an isomorphism of $\gxg$-modules. As a consequence,
$\Lambda^+_{\cO}$ identifies to $\Lambda'\cap\Lambda^+$.

\item The space $k[\cO]^{\diag G}$ of conjugation-invariant
functions has a basis consisting of the functions
$\chi_{\lambda,K}:(x,y)\mapsto\Tr_{V_\lambda}(xp_{\lambda,K}y^{-1})$,
where $\lambda\in\Lambda'\cap \Lambda^+$.

\item The multiplication in $k[\cO]$ is given by
$$
f_{\lambda}\cdot f_{\mu} = 
\sum_{\nu\in\Lambda^+,\,\nu\leq_K\lambda+\mu} 
\Tr_1(p^{\nu}_{\lambda\mu}(f_{\lambda}\otimes f_{\mu})),
$$
for any $f_\lambda \in \End V_\lambda$ and $f_\mu\in \End V_\mu$.
In particular,
$$
\chi_{\lambda,K} \cdot \chi_{\mu,K} = 
\sum_{\nu\in\Lambda^+,\,\nu\leq_K\lambda+\mu} 
c_{\lambda\mu}^{\nu} \chi_{\nu,K},
$$
and the product $\End V_\lambda \cdot \End V_\mu$ is the sum of those
$\End V_\nu$ such that $\nu\leq_K\lambda+\mu$ and
$c_{\lambda\mu}^{\nu}\neq 0$.

\end{enumerate}

\end{proposition}

\begin{proof}
(1) Clearly, $p_{\lambda,K}$ is fixed by 
$(R_u(P_K)\times R_u(P_K^-))\diag L_K$, and hence by 
$(R_u(P_I)\times R_u(P_I^-))\diag L_K$ (since $K\subseteq I)$. 
Moreover, since $\lambda$ is orthogonal to $J$, it extends to a
character of $L_J$, and $L_J$ acts on $V_{\lambda,K}$ via that
character; thus, $H_J$ fixes $V_{\lambda,K}$ pointwise. 
As another consequence, $L_J$ stabilizes $V_{\lambda,K}$ and hence
$V_{\lambda}^K$ as well. It follows that $p_{\lambda,K}$ is fixed by 
$H_J\times H_J$. We conclude that $p_{\lambda,K}$ is fixed by $H$.

(2) One has an equivariant isomorphism
$$
\bigoplus_{\lambda,\mu\in \Lambda^+}
V_\lambda^* \otimes_k V_\mu \otimes_k V_\lambda\otimes_k V_\mu^*
\to k[\gxg]
$$
that maps $u\otimes v$ to the function $(x,y)\mapsto ((x,y)v)(u)$ ;
here $u\in V_\lambda^* \otimes_k V_\mu$, 
$v\in V_\lambda\otimes_k V_\mu^*=(V_\lambda^* \otimes_k V_\mu)^*$, and
$(x,y) \in \gxg$. This yields an isomorphism of $\gxg$-modules:
$$
\bigoplus_{\lambda,\mu\in \Lambda^+}
V_\lambda^* \otimes_k V_\mu \otimes_k (V_\lambda\otimes_k V_\mu^*)^H
\to k[\gxg]^H=k[\cO].
$$
Moreover, $(V_\lambda\otimes_k V_\mu^*)^H$ is the subset
of 
$(H_J\times H_J)\diag L_K$-invariants in
$$
(V_\lambda\otimes_k V_\mu^*)^{R_u(P_I)\times R_u(P_I^-)}=
V_{\lambda}^{R_u(P_I)}\otimes_k (V_\mu^*)^{R_u(P_I^-)}.
$$
And $V_{\lambda}^{R_u(P_I)}=V_{\lambda,I}$, whereas 
$(V_{\mu}^*)^{R_u(P_I^-)}$ is the dual of $V_{\mu,I}$ (since both are
simple $L_I$-modules with lowest weight $-\mu$). Thus, we obtain
$$
(V_\lambda\otimes_k V_\mu^*)^H =
\Hom^{(H_J\times H_J)\diag L_K}
(V_{\mu,I},V_{\lambda,I}).
$$
By Schur's lemma, that space is one-dimensional if both $V_{\lambda,I}$
and $V_{\mu,I}$ are fixed pointwise by $H_J$, and are
isomorphic as (simple) $L_K$-modules. This is equivalent to:
$\lambda=\mu$ vanishes on $T'$, that is,
$\lambda=\mu\in\Lambda'$; then $\lambda$ is orthogonal to $J$, and 
$(V_{\lambda}\otimes_k V_{\lambda}^*)^H$
is spanned by $p_{\lambda,I}=p_{\lambda,K}$. Otherwise,
$(V_\lambda\otimes_k V_\mu^*)^H$ is zero. Now the assertion on the
structure of $k[\cO]$  follows by keeping track of the
isomorphisms.

(3) follows readily from (2).

(4) Let $f_\lambda$, $f_\mu \in k[\cO]$ be the images of
$u_\lambda\in \End V_\lambda$, $u_\mu\in \End V_\mu$, where 
$\lambda,\mu\in \Lambda'\cap \Lambda^+$. Then one has for $x,y\in G$:
$$\displaylines{
(f_\lambda\cdot f_\mu)(x,y)=
\Tr_{V_\lambda\otimes_k V_\mu}(u_\lambda x p_{\lambda,K}y^{-1}\otimes
u_\mu x p_{\mu,K}y^{-1})
\hfill\cr\hfill
=\Tr_{V_\lambda\otimes_k V_\mu}((u_\lambda\otimes u_\mu)(x\otimes x)
(p_{\lambda,K}\otimes p_{\mu,K})(y^{-1}\otimes y^{-1})).
\cr}$$
Now we claim that
$$
p_{\lambda,K}\otimes p_{\mu,K}=
\sum_{\nu\in\Lambda^+,\,\nu\leq_K \lambda+\mu}
1\otimes p_{\nu,K}
$$
in $\End V_\lambda\otimes_k \End V_\mu=\bigoplus_\nu
\End N_{\lambda\mu}^{\nu}\otimes_k \End V_\nu$. To see this, notice that
the subspace 
$$
V_{\lambda,K}\otimes_k V_{\mu,K}\subseteq V_\lambda\otimes_k V_\mu
$$ 
is the sum of the $T$-weight spaces associated with those weights
$\chi$ such that $\chi\leq_K\lambda+\mu$, so that 
$$
V_{\lambda,K}\otimes_k V_{\mu,K}=
\bigoplus_{\nu\in\Lambda^+,\,\nu\leq_K\lambda+\mu}
N_{\lambda\mu}^{\nu}\otimes_k V_{\nu,K}.
$$
The claim follows, since 
$$
p_{\lambda,K}\otimes p_{\mu,K}:V_\lambda\otimes_k V_\mu\to
V_{\lambda}\otimes_k V_{\mu}
$$
is the unique $T$-equivariant projection to 
$V_{\lambda,K}\otimes_k V_{\mu,K}$.

Using that claim, we obtain 
$$
(f_\lambda\cdot f_\mu)(x,y) = 
\sum_{\nu\in\Lambda^+,\,\nu\leq_K\lambda+\mu}
\Tr_{N_{\lambda\mu}^{\nu}\otimes_k V_\nu}
(p_{\lambda\mu}^{\nu}(u_\lambda\otimes u_\mu)
(1\otimes xp_{\nu,K}y^{-1}))
$$
which implies our assertions.
\end{proof}

Next we characterize those homogeneous spaces
$\gxg/H_{K,\Lambda',J}$ that are quasi-affine, and we describe their
$\diag T$-fixed points.

\begin{proposition}\label{quasiaffine}
Let $H=H_{K,\Lambda',J}$ and $\gxg/H=\cO$. Then $\cO$ is
quasi-affine if and only if: $\Lambda'\cap \Lambda^+$ spans
$\Lambda'$, and $J$ is the set of those simple roots that are
orthogonal to $\Lambda'$.

In that case, one has $\Lambda_{\cO}=\Lambda'$, and 
$\cO^{\diag T}=(\diag W)\cO^{\diag T}_0$,
where $\cO^{\diag T}_0$ is the $T\times T$-orbit 
$(T\times T)H/H\simeq T/T'$. Moreover, the subgroup of $W$ that fixes
$\cO^{\diag T}_0$ pointwise (resp.~stabilizes $\cO^{\diag T}_0$) is
$W_J$ (resp.~$W_J\times W_K$). 

Finally, $(\diag G)\cO^{\diag T}$ is dense in $\cO$.
\end{proposition}

\begin{proof}
Since the monoid $\Lambda'\cap\Lambda^+$ is finitely generated, the
algebra $k[\cO]$ is finitely generated as well (this follows
e.g. from \cite{Grosshans97} Theorem 16.2). Let $X$ be the corresponding
affine irreducible $\gxg$-variety, then one has a dominant morphism 
$\cO\to X$ with irreducible general fibers, since
$k[\cO]=k[\gxg]^H$ is integrally closed in the function field
$k(\cO)=k(\gxg)^H$. Thus, $\gxg/H$ is quasi-affine if and only
if $\dim \gxg/H=\dim X$.

Now $\dim X=\rank(\Lambda'\cap \Lambda^+) + 
\vert \Phi - (\Lambda'\cap \Lambda^+)^{\perp}\vert$,
where $(\Lambda'\cap \Lambda^+)^{\perp}$ denotes the set of those
roots that are orthogonal to $\Lambda'\cap \Lambda^+$ (as follows
e.g. from Weyl's dimension formula). Then 
$(\Lambda'\cap \Lambda^+)^{\perp}=\Phi_{J'}$ for some
$J'\supseteq J$. So $\dim \cO = \dim X$ if and only if:
$\Lambda'\cap \Lambda^+$ spans $\Lambda'$, and $J'=J$.
This proves the first assertion.

For the second assertion, let $(x,y)\in\gxg$ such that the coset
$(x,y)H$ is in $\cO^{\diag T}$. We claim that 
$(x,y)\in \diag N (T\times T)H$. Note that  
$$
(nt_1 x h_1,nt_2 y h_2)\in\cO^{\diag T}
$$ 
for all $n\in N$, $(t_1,t_2)\in T\times T$ and $(h_1,h_2)\in H$. 
Moreover, since $H\subseteq P_I\times P_I^-$, the coset
$(xP_I,yP_I^-)$ is a $\diag T$-fixed point in $G/P_I\times G/P_I^-$.
Thus, $x\in NP_I$ and $y\in N P_I^-$. Since $P_I^-$ is the image of
$H$ under the second projection to $G$, we may assume that $y=1$.
Moreover, replacing $x$ by $xh_1$ where 
$h_1\in H \cap (G\times 1)=R_u(P_I)H_J$, we may assume that 
$x\in NL_K$. Then $(x^{-1}tx,t)\in L_K\times T$ for all $t\in T$. 
On the other hand, $(x^{-1}tx,t)\in H$ since $\diag T$ fixes $(x,1)H$.
Thus, 
$$
(x^{-1}tx,t)\in H\cap (L_K\times T)=(T'\times T')\diag T.
$$
In particular, $x^{-1}tx\in T$ for all $t\in T$, that is, $x\in N$. 
And $x^{-1}txt^{-1}\in T'$, so that $\chi(x^{-1}tx)=\chi(t)$ for all
$\chi\in\Lambda'$: in other words, $x$ fixes $\Lambda'$
pointwise. Thus, $x$ fixes pointwise $\Lambda'\cap \Lambda^+$,
so that $x\in N\cap L_J$. Hence $x\in T(H\cap(G\times 1))$, so that
$(x,1)\in(T\times T)H$. This completes the proof of the claim. 

By that claim, one has 
$\cO^{\diag T}= (\diag W)\cO^{\diag T}_0$. Moreover, a given 
$w=nT\in W$ stabilizes $\cO^{\diag T}_0$ if and only if 
$(n,n)\in (T\times T)H$, that is, $w\in W_J\times W_K$. Finally, 
$w$ fixes $\cO^{\diag T}_0$ pointwise if and only if it fixes
$\Lambda'$ pointwise, that is, $w\in W_J$.

For the final assertion, note that 
$(\diag G)\cO^{\diag T}=(\diag G)\cO^{\diag T}_0$ is dense
in $\cO$ if and only if the subset $(\diag G) (T\times T) H$ is dense
in $\gxg$, or (quotienting by $\diag G$): the subset
$$
\{xty^{-1}~\vert~(x,y)\in H,~t\in T\}
$$ 
is dense in $G$. Since the product $R_u(P_I)L_I R_u(P_I^-)$ is open in
$G$, it suffices to show that the subset 
$$
\{xty^{-1}~\vert~(x,y)\in H\cap(L_I\times L_I),~t\in T\}
$$ 
is dense in $L_I$. Since $L_I=[L_J,L_J]L_K$ and 
$H\cap (L_I\times L_I)=(H_J\times H_J)\diag L_K$ where $H_J$ contains
$[L_J,L_J]$, this reduces to the well-known fact that the subset 
$\{hth^{-1}~\vert~h\in L_K,~t\in T\}$ is dense in $L_K$.
\end{proof}

Finally, we record the following easy result on the multiplication in
the ring $k[\gxg/H_{K,\Lambda',J}]$.

\begin{lemma}\label{transvectants}
Let $H=H_{K,\Lambda',J}$ be such that $\cO=\gxg/H$ is quasi-affine. 
Then, for any $\alpha\in K$, there exists
$\lambda\in\Lambda'\cap\Lambda^+$ such that
$\langle\lambda,\check\alpha\rangle\neq 0$. In that case,
$2\lambda-\alpha\in\Lambda^+$, and the product 
$\End V_\lambda\cdot\End V_\lambda\subset k[\cO]$ contains a unique
copy of the $\gxg$-module $\End V_{2\lambda-\alpha}$. A corresponding
$B^-\times B$-eigenvector is
$$
(\eta_\lambda\otimes v_\lambda)\cdot
(e_{\alpha}\eta_\lambda\otimes e_{-\alpha}v_\lambda)
-(e_{\alpha}\eta_\lambda\otimes v_\lambda)\cdot
(\eta_\lambda\otimes e_{-\alpha}v_\lambda),
$$
where $\eta_\lambda$ (resp.~$v_\lambda$) is a lowest (resp.~highest)
weight vector in $V_\lambda^*$ (resp.~$V_\lambda$), and
$e_{\pm\alpha}$ is a root vector in $\Lie(G)$ of weight $\pm\alpha$.
\end{lemma}

\begin{proof}
If $\check\alpha$ is orthogonal to $\Lambda'\cap\Lambda^+$, then it is 
orthogonal to $\Lambda'$ by Proposition \ref{quasiaffine}. But this
contradicts the assumption that $\Lambda'$ contains $K$. This proves
the first assertion.

For any $\lambda\in\Lambda^+$ and $\alpha\in\Pi$ such that
$\langle\lambda,\check\alpha\rangle\neq 0$, the weight
$2\lambda-\alpha$ is dominant, and the corresponding simple module
$V_{2\lambda-\alpha}$ occurs in $V_\lambda\otimes_k V_\lambda$ with
multiplicity $1$, a corresponding highest weight vector being
$v_\lambda\otimes e_{-\alpha}v_\lambda - 
e_{-\alpha}v_\lambda \otimes v_\lambda$.
This implies the remaining assertions.
\end{proof}


\section{Reductive varieties}

\label{sec: Reductive varieties}


\subsection{Classification}

\label{subsec: Classification1}

{}From the combinatorial classification of the group-like homogeneous
spaces in Section \ref{sec: structure}, we will deduce a similar
classification of reductive varieties. We begin with a
characterization of these varieties in terms of their weight set.

\begin{proposition}\label{irreducible}
For an affine irreducible $\gxg$-variety $X$, the following conditions 
are equivalent: 

\begin{enumerate}

\item $X$ is a reductive variety.

\item $X$ is diagonal, multiplicity-free, and $\Lambda_X^+$ is
  saturated in $\Lambda$. 
\end{enumerate}

Then all $\gxg$-orbit closures in $X$ are normal.
\end{proposition}

\begin{proof}
(1)$\Rightarrow$(2) By Proposition \ref{module}, $X$ is diagonal and
multiplicity-free. And since the isotropy groups of
$X$ are connected, $\Lambda_X$ is saturated in $\Lambda$
by Lemma \ref{automorphisms} and Proposition \ref{quasiaffine}. 
Together with Lemmas \ref{multfree} and \ref{normality}, it follows
that all $\gxg$-orbit closures in $X$ are normal, and that
$\Lambda_X^+=\Lambda\cap C_X$.

(2)$\Rightarrow$(1) follows from Propositions \ref{connected} and
\ref{module}.

\end{proof}

Next we construct all reductive varieties. For any subset 
$K\subseteq\Pi$, consider the closed subgroup of $\gxg$:
$$
H_K=H_{K,\Lambda,\emptyset}=(R_u(P_K)\times R_u(P_K^-))\diag L_K.
$$
Let $k[G]_{(K)}$ be the algebra of regular functions on the
homogeneous space $\gxg/H_K$. Then 
$k[G]_{(K)}=\bigoplus_{\lambda\in\Lambda^+} \End V_\lambda$ 
as a $\gxg$-module; its multiplication is given by
$$
f_{\lambda}\cdot f_{\mu} = 
\sum_{\nu\in\Lambda^+,\,\nu\leq_K\lambda+\mu} 
\Tr_1(p^{\nu}_{\lambda\mu}(f_{\lambda} \otimes f_{\mu})).
$$
Thus, we may regard the algebra $k[G]_{(K)}$ as a degeneration of
$k[G]$; this will be developed in Subsection 
\ref{subsec: The Vinberg family} below.  

\begin{proposition}\label{classification1}
Let $K$ be a subset of $\Pi$ and let $C$ be a cone in
$\Lambda_{\bR}^+$, satisfying the following conditions:

\begin{enumerate}

\item The linear span $\lin C$ contains $K$.

\item $C$ is the intersection of $\lin C\cap\Lambda^+_{\bR}$ with
a finite (possibly empty) family $(v_i\geq 0)_{i\in I}$ of closed
half-spaces of $\Lambda_{\bR}$, where $v_i\leq 0$ on $K$
for all $i$.

\end{enumerate}

Then $C$ is a rational polyhedral convex cone; the $\gxg$-invariant
subspace
$$
k[G]_{(C,K)}=\bigoplus_{\lambda\in\Lambda\cap C}\End V_{\lambda}
$$ 
of $k[G]_{(K)}$ is a finitely generated subalgebra, and the
corresponding affine $\gxg$-variety $X_{C,K}$ is a reductive variety
for $G$. It is equipped with a base point $x_{C,K}$ in its open 
$B^-\times B$-orbit, fixed by $\diag T$.

Conversely, each reductive variety is isomorphic to some
$X_{C,K}$, where $C,K$ are uniquely determined and satisfy (1) and
(2). Moreover, the isotropy group of $x_{C,K}$ is
$H_{K,\Lambda',J}$, where $\Lambda'$ is the subgroup
of $\Lambda$ generated by $\Lambda\cap C$, and $J$ is the set of
simple roots orthogonal to $C$. Finally, $\Aut^{\gxg}(X_{C,K})$ is a
diagonalizable group with character group $(\Lambda\cap\lin C)/\bZ K$.

\end{proposition}

\begin{proof}
This follows from embedding theory of spherical homogeneous spaces
(see \cite{Knop89}). We give a direct proof for the reader's
convenience. 

The multiplication in $k[G]_{(K)}$ satisfies 
$$
\End V_{\lambda}\cdot \End V_{\mu}\subseteq
\bigoplus_{\nu\leq_K \lambda+\mu} \End V_{\nu}.
$$
Together with condition (2), it follows that $k[G]_{(C,K)}$ is a
subalgebra of $k[G]_{(K)}$. Its subalgebra 
$k[G]_{(C,K)}^{U^-\times U}$ is isomorphic to the algebra of the
monoid $\Lambda\cap C$ over $k$. This monoid is 
finitely generated by Gordan's lemma; hence the same holds for
$k[G]_{(C,K)}^{U^-\times U}$, and for $k[G]_{(C,K)}$ by
\cite{Grosshans97} Theorem 16.2. Now $X_{C,K}$ is a reductive variety
by Proposition \ref{irreducible}. It is equipped with a dominant
$\gxg$-equivariant morphism $\gxg/H_K\to X_{C,K}$ and hence with a
base point $x_{C,K}$ (the image of the coset $H_K$).

Conversely, let $X$ be a reductive variety. Then
$X$ is normal by Proposition \ref{irreducible} again; and by
Propositions \ref{structure} and \ref{quasiaffine}, $X$ contains an
open $\gxg$-orbit isomorphic to $\gxg/H$ for some
$H=H_{K,\Lambda',J}$, where $J$ is the set of simple roots orthogonal
to $\Lambda'$. Then $H\subseteq H_K$ so that $k[X]$
identifies to a $\gxg$-invariant subalgebra of $k[G]_{(K)}$. Since
$X$ is normal, we have by Lemmas \ref{Uinv} and \ref{multfree}:
$k[X]=\bigoplus_{\lambda\in\Lambda\cap C}\End V_{\lambda}$
for a uniquely defined rational polyhedral convex cone $C$ in
$\Lambda_{\bR}^+$. Moreover, $\Lambda\cap C$ spans $\Lambda'$ (since
$k[X]$ and $k[\gxg/H]$ have the same fraction field), so that
$\lin C$ contains $K$. Let $(X_i)_{i\in I}$ be the (finite, possibly
empty) set of closed irreducible $\gxg$-subvarieties of codimension
$1$ in $X$, and let $(v_i)_{i\in I}$ be the corresponding set of
valuations of the field of rational functions $k(X)=k(\gxg/H)$. 
Then each valuation ring $\cO_{v_i}$ is stable by $\gxg$, and  
$$
k[X]=k[\gxg/H]\cap\bigcap_{i\in I}\cO_{v_i}.
$$ 
Moreover, for any $T\times T$-eigenvector $f\in k[X]^{U^-\times U}$
with weight $\lambda$, the value $v_i(f)$ depends only on $\lambda$
(since $\lambda$ determines uniquely $f$ up to scalar). This defines
an additive map $v_i:\Lambda^+_X\to \bZ$; it extends to a linear
functional on $\Lambda_{\bR}$ that we still denote by $v_i$. Then
$$
\Lambda_X^+=\Lambda_{\gxg/H}\cap \bigcap_{i\in I}(v_i\geq 0).
$$
Finally, since $v_i$ is a $\gxg$-invariant valuation, it is constant
on every $\End V_\lambda -\{0\}$, with value $v_i(\lambda)$. And since 
$$
\End V_{2\lambda-\alpha}\subset \End V_\lambda \cdot \End V_\lambda
$$
for any $\alpha\in \Pi$ such that 
$\langle\lambda,\check\alpha\rangle>0$ (Lemma \ref{transvectants}), 
it follows that $v_i(\alpha)\leq 0$ for any $\alpha\in K$.

Recall from Lemma \ref{automorphisms} that
$\Aut^{\gxg}(\gxg/H)=Z_K/T_H$ is diagonalizable with character group
$\Lambda'/\bZ K$. Moreover, $\Aut^{\gxg}(\gxg/H)$ acts by scalars on
every subspace $\End V_\lambda$ of $k[\gxg/H]$ (the corresponding
character being $\lambda\in\Lambda'$.) Thus, $\Aut^{\gxg}(\gxg/H)$
stabilizes the subalgebra $k[X]$; this yields an isomorphism 
$\Aut^{\gxg}(X) \simeq \Aut^{\gxg}(\gxg/H)$. 
\end{proof}


\subsection{Associated stable toric varieties}

\label{subsec: Associated stable toric varieties}

We will describe the subsets of $\diag T$-fixed points in reductive
varieties for $G$. A useful preliminary result is the following:

\begin{lemma}\label{closed}
Let $X$ be a stable reductive variety for $G$. Then $X$ contains a
unique closed $\gxg$-orbit $Y$. The latter is isomorphic to $G/G_1$,
where $G_1$ is a closed connected normal subgroup of $G$ (and $\gxg$
acts on $G/G_1$ by left and right multiplication). In particular,
$Y^{\diag T}$ is a unique $T\times T$-orbit.

Moreover, there exists a stable reductive variety $X_1$ (for $G_1$)
with a fixed point of $G_1\times G_1$, such that 
$X\simeq G\times^{G_1} X_1$ (the quotient of $G\times X_1$ by the
$G_1$-action: $g_1\cdot(g,x)=(gg_1^{-1}, g_1x)$.)
\end{lemma}

\begin{proof}
Since $X$ is connected and contains only finitely many $\gxg$-orbits,
$X//G$ is connected and finite. Thus, $X$ contains a unique closed
orbit $Y=\gxg/H$; then $H$ is connected and reductive. By Proposition
\ref{structure}, it follows that $H=(H_J\times H_J)\diag L_K$, where 
$\Pi$ is the union of orthogonal subsets $J$ and $K$, and 
$[L_J,L_J]\subseteq H_J\subseteq [L_J,L_J] Z_K$. Thus, $H_J$ is a
normal subgroup of $G$, and $H_J\times H_J$ fixes pointwise
$\gxg/H$. Replacing $G$ by $G/H_J$, we may assume that $H_J$ is
trivial; then $L_K=G$, so that $Y\simeq G$. This proves the first
assertion; the second assertion then follows from Luna's slice
theorem.
\end{proof}

\begin{proposition}\label{fixed}
Let $X=X_{C,K}$ be a reductive variety for $G$, with open 
$B^-\times B$-orbit $X_0$. Let $X'$ be the closure of $X^{\diag T}_0$
in $X$. Then:

\begin{enumerate}
\item  $X'$ is a toric variety for $T$, with lattice $\Lambda_X$
and cone $W_K C$. In particular, $X'$ is normal, and $\sigma=W_K C$ is
a rational polyhedral convex cone. Moreover, one has
$$
\sigma\cap\Lambda^+_{\bR}=C, \text{ and }
K=\{\alpha\in\Pi~\vert~
\sigma^0\text{ meets the hyperplane } (\alpha=0)\},
$$
where $\sigma^0$ is the relative interior of $\sigma$.

\item $X^{\diag T}=(\diag W) X'$ is a diagonal, multiplicity-free 
$T$-variety with weight set $\Lambda\cap W\sigma$. As a consequence, 
$X^{\diag T}$ is a stable toric variety, and the distinct translates
$w\sigma^0$ $(w\in W$) are disjoint. Moreover, $X^{\diag T}$ meets any
$\gxg$-orbit along a unique $WT$-orbit.

\end{enumerate}

\end{proposition}

\begin{proof}
(1) Let $x=x_{C,K}\in X_{C,K}$ be the base point, with isotropy group 
$H$; then $x\in X_0^{\diag T}$. Thus, $X'$ 
is the closure of $(T\times T)x$, so that we may identify  
$k[X']$ with a $T\times T$-invariant subalgebra of 
$k[T\times T/\diag T]\simeq k[T]$. On the other hand, the vector
space $k[X']$ is spanned by the restrictions to $X'$ of all
subspaces $\End V_{\lambda}$ of $k[X]$. By Proposition
\ref{product}, it follows that $k[X']$ is spanned by the functions
on $T\times T$:
$$
(x,y)\mapsto \Tr_{V_{\lambda}}(uxp_{\lambda,K}y^{-1})
=\Tr_{V_{\lambda,K}}(p_{\lambda,K}y^{-1}ux)
$$
where $u\in \End V_{\lambda}$.
So the weights of $\End V_{\lambda}\vert_{X'}$ are
precisely the weights of the $T$-module $V_{\lambda,K}$, that is,
those weights $\chi\in\Conv(W_K\lambda)$ that satisfy
$\chi\leq_K\lambda$. This set of weights is $W_K$-invariant, and for
any such weight $\chi$, there exists $w\in W_K$ such that 
$\langle w\chi,\check\alpha\rangle \geq 0$ for all $\alpha\in K$. 
On the other hand, since $\lambda$ is dominant and 
$w\chi\leq_K \lambda$, one has
$\langle w\chi,\check\alpha\rangle \geq 0$ for all 
$\alpha \in \Pi - K$. So $w\chi\in\Lambda^+_{\bR}$, whence 
$w\chi\in \Lambda\cap C$ by condition (2) of Proposition
\ref{classification1}. Thus, 
$\Lambda^+_{X'}\subseteq W_K(\Lambda\cap C)$, and the opposite 
inclusion is obvious: we have proved that 
$\Lambda^+_{X'} = W_K(\Lambda\cap C)$. 

It follows that $\Lambda_{X'}=\Lambda_X$ and that $\Lambda_{X'}^+$ is 
saturated in this subgroup. Since $X'$ is irreducible, it is a toric
variety with cone $W_K C=\sigma$. Moreover, the argument shows that
$\sigma\cap\Lambda^+_{\bR}=C$. 

Let $\alpha\in \Pi$. If $\alpha\notin K$, then $w\alpha\in\Phi^+$ for
all $w\in W_K$. Thus, $\alpha\geq 0$ on $W_K\Lambda^+_{\bR}$, and
hence on $\sigma$. Conversely, if $\alpha\in K$, then $\alpha$ takes
positive values on $C$ (otherwise $\alpha$ is orthogonal to $C$, which
contradicts the assumption that $\alpha\in \lin C$). Thus, $\alpha$
takes negative values on $\sigma=W_K C$, so that $C$ meets the
hyperplane $(\alpha=0)$. This completes the proof of (1). 

(2) By Proposition \ref{quasiaffine}, $X^{\diag T}$ meets any
$\gxg$-orbit $\cO$ along a unique $WT$-orbit. On the other hand, we
assert that $X'$ meets any $\gxg$-orbit. This may be deduced from
embedding theory of spherical homogeneous spaces (see \cite{Knop89});
here is a direct argument, using the notation of Proposition
\ref{structure}. Since $P_I\times P_I^-$ is a parabolic subgroup of
$\gxg$, the closure of $(P_I\times P_I^-)x$ in $X$ meets all
$\gxg$-orbits. Moreover, 
$$
(P_I\times P_I^-)x \simeq (P_I\times P_I^-)H/H\simeq 
(L_I\times L_I)/(L_I\times L_I)\cap H
$$
is isomorphic (as a $L_I\times L_I$-variety) to the quotient
$L_I/H_J$ where $L_I$ acts by left and right multiplication. This 
quotient is a connected reductive group $G_1$, 
and the image of $T$ is a maximal torus $T_1$. Now the closure of $T_1$
in any $G_1\times G_1$-equivariant embedding of $G_1$ meets all
$G_1\times G_1$-orbits, as follows e.g. from a theorem of Iwahori
(see \cite{Mumford_GIT3ed}, p.~52).

It follows that the irreducible components of $X^{\diag T}$ are the
$wX'$ ($w\in W$), and also (using Lemma \ref{closed}) that 
$X^{\diag T}$ is connected. Thus, the weight set of $X^{\diag T}$ is
$W(\Lambda\cap C)=\Lambda\cap W\sigma$. We now claim that 
$X^{\diag T}$ is multiplicity-free. It suffices to check that the
multiplicity of any $\lambda\in \Lambda\cap C$ is $1$. Let $f$ be a 
$T\times T$-eigenvector in $k[X]^{U^-\times U}$ of weight
$\lambda$, then $f\vert_{X'}$ is a $T\times T$-eigenvector of weight
$\lambda$, which does not vanish identically on any $wX'$, 
$w\in \Stab_W(\lambda)$. Let $\varphi\in k[X']$ be another 
$T\times T$-eigenvector of weight $\lambda$. If $\varphi$ does not
vanish identically on $wX'$, then $w^{-1}(\lambda)\in W_K C$, whence
$w\in \Stab_W(\lambda) W_K$, and $wX'\subseteq \Stab_W(\lambda) X'$. 
So it suffices to show that $\varphi$ is a scalar multiple of $f$ when 
restricted to $\bigcup_{w\in \Stab_W(\lambda)} wX_0^{\diag T}$, a
subset of $X_f^{\diag T}$. With the notation of Lemma \ref{local},  
$X_f^{\diag T}=Z^{\diag T}$; the latter is connected and contains
only finitely many $T\times T$-orbits. Thus, the quotient
$\varphi/f$, a regular $T\times T$-invariant function on
$Z^{\diag T}$, is constant. This proves our claim.

Now $X^{\diag T}$ is a stable toric variety, by the claim and 
Lemma \ref{saturated}. And since it is multiplicity-free, the
distinct $w\sigma^0$ are disjoint.
\end{proof}

\begin{proposition}\label{invariants}
For any reductive variety $X$, the natural morphism 
$$
p:X^{\diag T}/\diag W\to X//\diag G
$$ 
is an isomorphism. 

Moreover, the closed $\diag G$-orbits in $X$ are
exactly the orbits of points in $X^{\diag T}$; the union of these
orbits contains a dense subset of every $\gxg$-orbit.

Finally, $X$ is self-adjoint for a unique automorphism $\Theta_X$
fixing pointwise $X^{\diag T}$.
\end{proposition}

\begin{proof}
By the main result of \cite{Luna75}, the morphism $p$ is finite; on
the other hand, it is well-known that the orbit $(\diag G)x$ is closed
for every $x\in X^{\diag T}$, and that its $T$-fixed point subset is
just $(\diag W)x$. Moreover, by Proposition \ref{quasiaffine}, the set 
$(\diag G) X^{\diag T}$ contains a dense open subset of every
$\gxg$-orbit. As a consequence, $p$ is birational. Since $X$ is
normal, $p$ is an isomorphism. In particular, it is surjective, so
that every closed $\diag G$-orbit meets $X^{\diag T}$.

By Lemma \ref{automorphisms}, the open orbit $\gxg/H$ is self-adjoint
for an automorphism $\Theta_{\gxg/H}$ fixing pointwise 
$(T\times T)H/H$. Since $\Theta(\End V_\lambda)\simeq\End V_\lambda$
as $\gxg$-modules, $\Theta_{\gxg/H}$ stabilizes each simple submodule
of $k[\gxg/H]$. Hence it stabilizes the subalgebra $k[X]$; this 
yields an automorphism $\Theta_X$ making $X$ self-adjoint. Since
$(\diag W)(T\times T)H/H$ is dense in $X^{\diag T}$, and $\Theta$
fixes $\diag W$ pointwise, it follows that $\Theta_X$ fixes pointwise
$X^{\diag T}$. Thus,
$$
\Theta_X((g,g)x)=(\theta(g),\theta(g))x
$$
for any $g\in G$ and $x\in X^{\diag T}$. This implies uniqueness of
$\Theta_X$, since $(\diag G)X^{\diag T}$ is dense in $X$.
\end{proof}

Next we reformulate the combinatorial classification of
reductive varieties by pairs $(C,K)$, in terms of the unique datum
$W_K C=\sigma$. This makes sense in view of Proposition \ref{fixed},
which motivates the following

\begin{definition}\label{admissible1}
A \emph{$W$-admissible cone} is a
rational polyhedral convex cone $\sigma$ in $\Lambda_{\bR}$,
satisfying the following conditions:

\begin{enumerate}

\item The relative interior $\sigma^0$ meets $\Lambda^+_{\bR}$.

\item The distinct $w\sigma^0$ ($w\in W$) are disjoint. 

\end{enumerate}

\end{definition}

Thus, if $C,K$ satisfy the conditions of Proposition
\ref{classification1}, then $\sigma=W_K C$ is a $W$-admissible
cone. Conversely, we have the following easy result.

\begin{lemma}\label{walls}
Let $\sigma$ be a $W$-admissible cone. Let
$C=C_\sigma=\sigma\cap\Lambda^+_{\bR}$ and 
$K=K_\sigma=\{\alpha\in\Pi~\vert~(\alpha=0) \text{ meets }\sigma^0\}$.
Then $C$, $K$ satisfy conditions (1), (2) of Proposition
\ref{classification1}.
\end{lemma}

\begin{proof}
Let $\alpha\in K$. Then $s_\alpha\sigma^0$ meets $\sigma^0$. Thus,
$s_\alpha\sigma^0=\sigma^0$, and $\sigma$ is invariant under
$s_\alpha$. So $\sigma$ is invariant under $W_K$, and also
$K\subset\lin \sigma = \lin C$. Hence condition (1) holds.
Moreover, $\sigma$ is the intersection of finitely many closed
half-spaces $(v_i\geq 0)_{i\in I}$, together with all their
$W_K$-translates. Moving $v_i$ in its $W_K$-orbit, we may
assume that $v_i\leq 0$ on $K$ for any $i\in I$.

We claim that $\sigma = W_K C$.  
Let $\tau=\sigma\cap(\alpha\geq 0)$, then 
$\sigma=\tau\cup s_\alpha\tau$, and 
$\tau^0\cap s_\alpha\tau^0=\emptyset$. It follows that $\tau$ is a
$W$-admissible cone, such that $C_\tau=C_\sigma$ and
$K_\tau \subseteq K_\sigma-\{\alpha\}$. By induction, we thus obtain
$\sigma\subseteq W_K C_\sigma$. The opposite inclusion is obvious;
this proves the claim. 

This claim implies that 
$C=\lin C\cap \Lambda^+_{\bR}\cap \bigcap_{i\in I}(v_i\geq 0)$; hence
condition (2) holds.
\end{proof}

We may now reparametrize the reductive varieties for $G$ in terms of
$W$-admissible cones, by setting
$$
X_{\sigma}=X_{C_\sigma,K_\sigma}.
$$
Then $X_\sigma^{\diag T}$ is the stable toric variety associated with
the complex of cones in $\Lambda_{\bR}$ consisting of all translates
$w\sigma$ and their faces; these are in bijection with the 
$T$-orbits in $X_\sigma^{\diag T}$ (see \cite{Alexeev_CMAV} Corollary
2.3.9). Thus, the $WT$-orbits in $X_\sigma^{\diag T}$ are in bijection
with the $W$-orbits of faces of $w\sigma$, $w\in W$. Moreover, the
restriction map
$$
\Aut^{\gxg}(X_\sigma) \rightarrow \Aut^{WT}(X_\sigma^{\diag T})
$$ 
is an isomorphism. 

Indeed, let $X=X_\sigma$, $K=K_\sigma$ and let $X'$ be as in
Proposition \ref{fixed}. Then $X^{\diag T}=WX'$, and 
$\Stab_W(X')=\Stab_W(\sigma)$ acts on it through its quotient $W_K$,
so that $\Aut^{WT}(X^{\diag T})$ injects into $\Aut^{W_K T}(X')$. 
Moreover, $\Aut^T(X')$ is the quotient of $T$ with character group 
$\Lambda\cap\lin \sigma$, and $\Aut^{W_K T}(X')$ is its subgroup with
character group $\Lambda\cap\lin \sigma/\bZ K$. By Proposition
\ref{classification1}, it follows that 
$\Aut^{\gxg}(X_\sigma)\to \Aut^{W_K T}(X')$
is an isomorphism.

So we have proved the following ``toric correspondence'':

\begin{theorem}\label{correspondence1}
The assignment $X\mapsto X^{\diag T}$ defines a bijection
from the reductive varieties for $G$, to the stable toric 
varieties for $T$ with a compatible $W$-action such that the quotient
by $W$ is irreducible. Moreover, the $\gxg$-orbits in $X$ are in
bijection with the $WT$-orbits in $X^{\diag T}$, and 
$\Aut^{\gxg}(X)$ is isomorphic to the automorphism group of the
$WT$-variety $X^{\diag T}$.
\end{theorem}


\section{Stable reductive varieties}

\label{sec: Stable reductive varieties}

We show how to obtain all stable reductive varieties, by glueing 
reductive varieties along invariant closed subvarieties. For this, we
introduce a combinatorial object, that will encode part of the glueing
data.

\begin{definition}\label{complex}
A \emph{$W$-complex of cones $\Sigma$ referenced by $\Lambda$} is a
topological space $\vert\Sigma\vert$ represented as a finite union of
distinct closed subsets $\sigma$ ($\sigma\in\Sigma$), together with a
map $\rho:\vert\Sigma\vert\to\Lambda_{\bR}$ such that:

\begin{enumerate}

\item $\rho$ identifies each $\sigma\in\Sigma$ with a rational
polyhedral convex cone in $\Lambda_{\bR}$.

\item If $\sigma\in\Sigma$, then each face $\tau\prec\sigma$ is in
$\Sigma$. 

\item If $\sigma$, $\tau$ in $\Sigma$, then their intersection in
$\vert\Sigma\vert$ is a union of faces of both.

\item $W$ acts on $\vert\Sigma\vert$, the reference map $\rho$ is
$W$-equivariant, and its restriction to any subset
$\bigcup_{w\in W} w\sigma$ is injective.

\end{enumerate}

(In particular, $W$ permutes the subsets $\sigma$.)
\end{definition}

For example, any $W$-admissible cone $\sigma$ defines a $W$-complex
of cones $\Sigma$, as follows: the cones in $\Sigma$ are the
$W$-translates of faces of $\sigma$, and the map
$\rho:\vert\Sigma\vert\to \Lambda_{\bR}$ is just the inclusion 
$\bigcup_{w\in W} w\sigma\subseteq \Lambda_{\bR}$. Such a $W$-complex
of cones will be called \emph{elementary}. Clearly, a $W$-complex of
cones is elementary if and only if it admits a unique maximal cone
modulo $W$.

Now consider an arbitrary $W$-complex of cones $\Sigma$, and the 
orbit space $\Sigma/W$. Note that any $\overline{\sigma}\in\Sigma/W$
admits a unique representative $\sigma\in\Sigma$ such that
$\rho(\sigma)^0$ meets $\Lambda^+_{\bR}$. This yields a partial
ordering $\leq$ on $\Sigma/W$, where 
$\overline{\tau}=W\tau\leq W\sigma=\overline{\sigma}$ if and only if 
$\tau \prec \sigma$. 

Every $\overline{\sigma}\in\Sigma/W$ defines a reductive variety
$X_{\overline{\sigma}}= X_{\sigma}$. Moreover,
$\overline{\tau}\leq\overline{\sigma}$ if and only if 
$\tau\prec\sigma$. In that case, by the results of Subsection
\ref{subsec: Associated stable toric varieties}, we have a closed
$\gxg$-equivariant immersion
$$
i_{\tau\sigma}: X_{\tau}\rightarrow X_{\sigma}
$$
such that the comorphism
$$
i_{\tau\sigma}^{\#}:
k[X_\sigma] = \bigoplus_{\lambda\in\Lambda^+\cap\sigma} \End V_\lambda
\rightarrow
\bigoplus_{\lambda\in\Lambda^+\cap\tau} \End V_\lambda = k[X_\tau]
$$
is the obvious projection.

Clearly, the $i_{\tau\sigma}$ define a directed system of
reductive varieties, indexed by the poset $\Sigma/W$.
We introduce additional twists of this system, as follows. Consider a
collection
$$
t=\{t_{\tau\sigma}\in \Aut^{\gxg}(X_{\tau}) ~\vert~ \tau\prec\sigma\}
$$ 
satisfying 
$$
t_{\sigma''\sigma}=t_{\sigma''\sigma'}\circ 
t_{\sigma'\sigma}|_{X_{\sigma''}}
$$ 
for all triples $\sigma''\prec \sigma' \prec \sigma$ (note that every
automorphism of $X_{\sigma'}$ leaves $X_{\sigma''}$ invariant). This
gives us a twisted directed system of reductive varieties
$$
t_{\tau\sigma} i_{\tau\sigma}:X_{\tau}\to X_{\sigma}.
$$

\begin{definition}
$X_{\Sigma,t} =\varinjlim X_{\sigma}$.
\end{definition}

Note that the twists are just the $1$-cocycles of the complex of
diagonalizable groups 
$$
0\to \bigoplus_{\overline{\sigma}\in\Sigma/W} \Aut^{\gxg}(X_{\sigma}) 
\to \bigoplus_{\overline{\tau}\prec\overline{\sigma}} 
\Aut^{\gxg}(X_{\tau}) \to \cdots
$$
with the obvious differential. We denote this complex by
$C^*(\Sigma/W,\Aut)$, with cocycle groups $Z^i(\Sigma/W,\Aut)$ and
cohomology groups $H^i(\Sigma/W,\Aut)$. As in the toric case (see
\cite{Alexeev_CMAV} Section 2.3), the following is easy to prove: 

\begin{proposition}\label{classification2}

\begin{enumerate}

\item $X_{\Sigma,t}$ is a stable reductive variety (for $G$). Its
irreducible components are the varieties $X_{\sigma}$ where
$\sigma\in\Sigma$ is a maximal cone such that $\sigma^0$ meets
$\Lambda^+_{\bR}$.

\item The $\gxg$-orbits in $X_{\Sigma,t}$ are in bijection with the
$W$-orbits of cones in $\Sigma$.

\item $\Aut^{\gxg}(X_{\Sigma,t}) = H^0(\Sigma/W,\Aut)$.

\item The set of isomorphism classes of $\gxg$-varieties
$X_{\Sigma,t}$ for a fixed $W$-complex $\Sigma$ is
$H^1(\Sigma/W,\Aut)$.

\end{enumerate}

\end{proposition}

Together with the results of Subsection 
\ref{subsec: Associated stable toric varieties}, this implies the
general version of the ``toric correspondence'':

\begin{theorem}\label{correspondence2}
The stable reductive varieties for $G$ are precisely the varieties
$X_{\Sigma,t}$, where $\Sigma$ is a $W$-complex of cones, and 
$t\in Z^1(\Sigma/W,\Aut)$. 

Thus, the assignment $X\mapsto X^{\diag T}$ defines a bijective
correspondence from the stable reductive varieties (for $G$), to the
stable toric varieties (for $T$) with a compatible $W$-action. This
correspondence preserves orbits and automorphism groups. 

Moreover, the closed $\diag G$-orbits in $X$ are exactly the orbits
in $X^{\diag T}$; the union of these orbits contains a dense subset of
every $\gxg$-orbit, and the natural morphism 
$X^{\diag T}/W\to X//\diag G$ is an isomorphism.

Finally, $X$ is self-adjoint for a unique automorphism fixing
pointwise $X^{\diag T}$.
\end{theorem}
\begin{proof}
  Let $(X_i=\overline{\cO}_i)$ be the finite set of 
  $\gxg$-orbit closures in a stable reductive variety $X$. For every $i$,
  $k[X_i]$ embeds into $k[\cO_i]$. Hence, $X_i$ is diagonal and
  multiplicity-free.  The connectedness of stabilizers implies that
  the lattice $\Lambda_{X_i}$ is saturated in $\Lambda$ (by Lemma 
  \ref{automorphisms}), and that the normalization morphism 
  $\nu_i:\tilde X_i\to X_i$ is bijective (by Lemma \ref{normality}). 
  Therefore, the normalizations are reductive varieties 
  $\tilde X_i =  X_{\sigma_i}$ and the glueing data $t$ for $X_i$ are
  the same as the glueing data for $\tilde X_i$. The inclusions
  between $X_i$'s  define on the set $\{\sigma_i\}$ the structure of
  a complex of cones $\Sigma$. The collection of morphisms 
  $\tilde X_i \to X$ defines a finite bijective morphism 
  $\pi:X_{\Sigma,t} \to X$.  Since $X$ is seminormal, $\pi$ is an
  isomorphism.
\end{proof}


\section{Self-adjoint stable reductive semigroups}

\label{sec: Self-adjoint stable reductive semigroups}

We study the relations between stable reductive varieties,
self-adjoint stable reductive semigroups, and reductive monoids. 
We begin with the following observation.

\begin{lemma}\label{monoid}
Any normal reductive monoid with unit group a quotient of $G$ by a
connected normal subgroup is a self-adjoint reductive semigroup for
$G$.

Conversely, any reductive semigroup for $G$ having an identity
element is a normal reductive monoid, with unit group a quotient of
$G$ by a connected normal subgroup.
\end{lemma}

\begin{proof}
Let $X$ be a normal reductive monoid with unit group $G/G_1$, where
$G_1$ is a connected normal subgroup. Then the $\gxg$-action on $X$
factors through an action of $G/G_1\times G/G_1$. Thus, we may replace
$G$ by $G/G_1$, and assume that the unit group is $G$.

By \cite{Rittatore98} Theorem 2, $X$ is affine. And since $G$ is dense
in $X$, the algebra $k[X]$ is a normal $\gxg$-invariant subalgebra of 
$k[G]$. Thus, $X$ is multiplicity-free and diagonal, and
$\Lambda_X^+$ is saturated in $\Lambda$. By Proposition
\ref{connected}, $X$ satisfies conditions (1) and (2) of Definition
\ref{stablereductive}. 

Let $m:X\times X\to X$ be the multiplication, with identity element
$1$. Clearly, $m$ satisfies the equivariance condition of Definition
\ref{stablesemigroup}. Thus, it factors through a morphism
$m//G:(X\times X)//G\to X$, the quotient being for the $G$-action
defined by $g\cdot(x_1,x_2)=((1,g)x_1,(g,1)x_2)$. Note that $m//G$ is
surjective (since $m(x,1)=x$) and birational (since 
$m^{-1}(1)=\{(x,y)\in G\times G~\vert~xy=1\}$ is a unique
$G$-orbit). Since $X$ is affine and normal, it follows that $m//G$ is
an isomorphism. 

Thus, $X$ is a reductive semigroup for $G$. Moreover, the involution
$\Theta_G:g\mapsto \theta(g^{-1})$ extends to an involution
$\Theta_X$ of $X$, by (the proof of) Proposition \ref{invariants}. 
Clearly, $\Theta_X$ makes $X$ self-adjoint.

Conversely, let $X$ be a reductive semigroup for $G$, equipped with an
identity element $1$. Then $(1g^{-1})(g1)=1$ for any $g\in G$, by
the equivariance condition. Thus, $g1$ has a left inverse. By
\cite{Rittatore98} Corollary 1, it follows that $g1$ is a unit, with
inverse $1g^{-1}$. As a consequence, $g1g^{-1}=1$, so that $g1=1g$ for
all $g\in G$. This implies that the map $g\mapsto g1$ is an algebraic 
group homomorphism from $G$ to the unit group $G(X)$. Moreover, $\gxg$
has only finitely many orbits in $X$ (by Proposition \ref{saturated}),
and hence in $G(X)$ (where $\gxg$ acts by left and right
multiplication). Therefore, the image of $G$ in $G(X)$ is a subgroup
of finite index. But $G(X)$ is a connected algebraic group, by
\cite{Rittatore98} Theorem 1. So $G(X)$ is the quotient of $G$ by a
normal subgroup, which must be connected since all isotropy groups of
$\gxg$ are.
\end{proof}

\begin{lemma}\label{square}
Let $X$ be a self-adjoint stable reductive semigroup for $G$. Then
$X$ is a multiplicity-free stable reductive variety for
$G$. Moreover, any closed $\gxg$-subvariety $Y$ is an ideal (that is,
$Y$ contains all products $xy$, $yx$ where $x\in X$, $y\in Y$) and is
invariant under the automorphism $\Theta_X$.
\end{lemma}

\begin{proof}
We have an equivariant isomorphism 
$k[X]\simeq (k[X]\otimes_k k[X])^G$.
It identifies the subalgebra $k[X]^{\gxg}$ to 
$(k[X]\otimes_k k[X])^{G\times G\times G}$ (where $G\times G\times G$ 
acts by $(g_1,g_2,g_3)\mapsto (g_1,g_2)\otimes(g_2,g_3)$.) The latter
contains the subalgebra $k[X]^{\gxg}\otimes_k k[X]^{\gxg}$, so
that
$$
k[X]^{\gxg}\otimes_k k[X]^{\gxg}\subseteq k[X]^{\gxg}
$$
as a subalgebra. Since $k[X]^{\gxg}$ is a finitely generated
$k$-algebra, this is only possible if $k[X]^{\gxg}=k$. Thus, all
multiplicities of the $\gxg$-module $k[X]$ are finite.

Write 
$$
k[X]=\bigoplus_{\lambda,\mu\in\Lambda} 
m_{\lambda,\mu} V_{\lambda}^*\otimes_k V_{\mu}
$$
where the $m_{\lambda,\mu}$ are non-negative integers. Then 
$m_{\lambda,\mu}=m_{\mu,\lambda}$ (since $X$ is self-adjoint.) Moreover, 
$$
m_{\lambda,\mu}=\sum_{\nu\in\Lambda^+} m_{\lambda,\nu} m_{\nu,\mu}
$$
since $X$ is isomorphic to $(X\times X)//G$. It follows that
$$
m_{\lambda,\lambda}=m_{\lambda,\lambda}^2+
\sum_{\nu\neq\lambda} m_{\lambda,\nu}^2
$$
so that $m_{\lambda,\lambda}\leq 1$ and $m_{\lambda,\nu}=$ for any
$\nu\neq\lambda$. Thus, $X$ is diagonal and multiplicity-free.
By Proposition \ref{module}, it follows that $X$ is a stable reductive
variety.

Write $k[X]=\bigoplus_{\lambda\in\Lambda^+_X} \End V_\lambda$.
By the equivariance condition, the comorphism $m^{\#}$ of the
multiplication $m$ sends every $\End V_\lambda$ to 
$\End V_\lambda\otimes_k\End V_\lambda$. Since $Y$ is
$\gxg$-invariant, its ideal $I_Y\subset k[X]$ is a partial sum of
simple submodules $\End V_\lambda$. It follows that 
$$
m^{\#}(I_Y)\subseteq I_Y\otimes_k I_Y = 
(k[X]\otimes_k I_Y)\cap (I_Y\otimes_k k[X]).
$$
In other words, $m((X\times Y)\cup(Y\times X)\subseteq Y$.
Likewise, since each $\End V_\lambda$ is invariant under $\Theta_X$,
the same holds for $I_Y$, and for $Y$.
\end{proof}

In particular, any self-adjoint reductive semigroup is a reductive
variety. We now establish the converse:

\begin{theorem}\label{semigroup}
Any reductive variety $X$ admits a structure of self-adjoint
reductive semigroup for $G$; then $X^{\diag T}$ is a subsemigroup,
which makes it a stable reductive semigroup for $T$. Moreover, 
$\Aut^{\gxg}(X)$ acts simply transitively on the set of all reductive
semigroup structures on $X$. 
\end{theorem}

\begin{proof}
Let $X=X_\sigma$ be a reductive variety and let $K=K_\sigma$. Write 
$$
k[X]=\bigoplus_{\lambda\in\Lambda^+\cap\sigma} \End V_\lambda
$$
and recall that the multiplication of $k[X]$ is defined by
$$
f_\lambda\cdot f_\mu=
\sum_{\nu\in\Lambda^+\cap\sigma,\,\nu\leq_K\lambda+\mu}
\Tr_1(p^{\nu}_{\lambda\mu}(f_{\lambda}\otimes f_{\mu})),
$$
for any $f_\lambda\in \End V_\lambda$, $f_\mu\in\End V_\mu$.
Let $m:X\times X\to X$ be a morphism satisfying the equivariance
condition of Definition \ref{stablesemigroup}. Then, as noted in the
proof of Lemma \ref{square}, the comorphism 
$m^{\#}:k[X]\to k[X]\otimes_k k[X]$ sends any $\End V_\lambda$
to $\End V_\lambda\otimes_k\End V_\lambda$. Let
$$
m^{\#}_{\lambda}:\End V_\lambda\to 
\End V_\lambda\otimes_k\End V_\lambda
$$
be the restriction; identify $\End V_\lambda$ to 
$V_\lambda^*\otimes_k V_\lambda$, and  
$\End V_\lambda\otimes_k\End V_\lambda$ to 
$V_\lambda^*\otimes_k V_\lambda\otimes_k V_\lambda^*\otimes_k V_\lambda$.
Then by equivariance, there exists $c_\lambda\in k$ such that
$$
m^{\#}_{\lambda}(\eta\otimes v)=c_\lambda \sum_i 
\eta\otimes v_i\otimes \eta_i\otimes v,
$$
where $(v_i)$ is a basis of $V_\lambda$ consisting of
$T$-eigenvectors, and $(\eta_i)$ is the dual basis. Moreover,
$c_\lambda\neq 0$ by the invariance condition.

We claim that $m^{\#}:k[X]\to k[X]\otimes_k k[X]$ is compatible
with the multiplication if and only if: 
$c_\lambda c_\mu=c_\nu$ for all $\lambda,\mu,\nu$ in
$\Lambda^+\cap\sigma$ such that the product 
$\End V_\lambda\cdot\End V_\mu$ contains $\End V_\nu$, that is:
$V_\lambda\otimes_k V_\mu$ contains 
$V_\nu$, and $\nu\leq_K\lambda+\mu$. To check this, we identify every
$\End V_\lambda$ to its dual space, via the bilinear form 
$(u,v)\mapsto \Tr(uv)$. This identifies $m^{\#}_\lambda$ to the 
dual of the map
$$
\End V_\lambda \otimes_k \End V_\lambda \to \End V_\lambda,~
\varphi\otimes\psi\mapsto c_\lambda \varphi\circ\psi.
$$
Now the compatibility of $m^{\#}$ with multiplication is equivalent to
the commutativity of the following diagram:
$$\CD
\End V_\lambda\otimes_k \End V_\mu 
@>{m_\lambda^{\#}\otimes m_\mu^{\#}}>>
\End V_\lambda \otimes_k \End V_\lambda \otimes_k \End V_\mu
\otimes_k \End V_\mu
\\
@V{p_{\lambda\mu}^{\nu}}VV 
@V{p_{\lambda\mu}^{\nu}\otimes p_{\lambda\mu}^{\nu}}VV\\
\End V_\nu 
@>{m_\nu^{\#}}>> 
\End V_\nu\otimes_k\End V_\nu. \\
\endCD
$$
This amounts to the commutativity of the dual diagram
$$
\CD
\End V_\nu\otimes_k\End V_\nu @>>> \End V_\nu \\
@VVV @VVV \\
\End V_\lambda\otimes_k\End V_\lambda\otimes_k\End V_\mu\otimes_k\End V_\mu
@>>> \End V_\lambda\otimes_k\End V_\mu. \\
\endCD
$$
For this, let $\varphi,\psi\in\End V_\nu$. Then 
$$
(p^{\nu}_{\lambda\mu})^*(\varphi)=1\otimes\varphi\in 
\End N_{\lambda\mu}^{\nu}\otimes_k\End V_\nu\subseteq 
\End V_\lambda\otimes_k \End V_\mu,
$$
so that 
$$
c_\lambda (p^{\nu}_{\lambda\mu})^*(\varphi)\circ
c_\mu (p^{\nu}_{\lambda\mu})^*(\psi) = c_\lambda c_\mu 1\otimes
(\varphi\circ\psi),
$$
whereas
$$
c_\nu (p^{\nu}_{\lambda\mu})^*(\varphi\circ\psi)=c_\nu 1\otimes
(\varphi\circ\psi).
$$
This completes the proof of the claim.

Now choose $c_\lambda=1$ for any $\lambda\in\Lambda^+\cap\sigma$. Then 
$m^{\#}$ is compatible with multiplication, by that claim. And 
one checks similarly that $m^{\#}$ is compatible with
the automorphism $\Theta_X^{\#}$ of $k[X]$, associated with $\Theta_X$.

Thus, we obtain a morphism $m:X\times X\to X$ satisfying the
equivariance condition. Since every 
$m_{\lambda}^{\#}$ defines an isomorphism 
$\End V_\lambda\simeq (\End V_\lambda\otimes_k \End V_\lambda)^G$,
the induced map $(X\times X)//G\to X$ is an isomorphism. 
Moreover, $m$ is associative, as follows from the definition of
$m^{\#}$ and associativity of the composition of endomorphisms.

By equivariance, $m:X\times X\to X$ maps 
$X^{\diag T}\times X^{\diag T}$ to $X^{\diag T}$. Thus, $X^{\diag T}$
is a subsemigroup of $X$. Clearly, it satisfies the equivariance
condition. To check the remaining invariance condition, it suffices to
show that the map 
$$
m^{\diag T,\#}_\lambda:k[X^{\diag T}]_\lambda\to
k[X^{\diag T}]_\lambda \otimes_k k[X^{\diag T}]_\lambda
$$
(induced by multiplication) is non-zero for any 
$\lambda\in W(\Lambda^+\cap\sigma)$, where 
$k[X^{\diag T}]_\lambda$ denotes the $\lambda$-weight space. By
$W$-equivariance, we may assume that
$\lambda\in\Lambda^+\cap\sigma$. Let $v_\lambda\in V_\lambda$ be a
highest weight vector (of weight $\lambda$), and 
$\eta_\lambda\in V_\lambda^*$ be a lowest weight vector (of weight
$-\lambda$). Then 
$\eta_\lambda\otimes v_\lambda\in \End V_\lambda\subset k[X]$
restricts to a non-zero function 
$\varphi_\lambda\in k[X^{\diag T}]_\lambda$.  
Since  
$$
m^{\#}(\eta_\lambda\otimes v_\lambda)=\sum_i 
\eta_\lambda\otimes v_i\otimes \eta_i\otimes v_\lambda,
$$
and $v_\lambda$ (resp.~$\eta_\lambda$) is the unique vector of weight
$\lambda$ (resp.~$-\lambda$) in the basis $(v_i)$ (resp.~$\eta_i)$, 
it follows that
$$
m^{\diag T,\#}_\lambda(\varphi_\lambda)=
\varphi_\lambda\otimes\varphi_\lambda.
$$

Finally, note that $\Aut^{\gxg}(X)$ preserves the decomposition of the
$\gxg$-module $k[X]$; thus, we may identify $\Aut^{\gxg}(X)$ to the
multiplicative group of families 
$(c_\lambda)_{\lambda\in\Lambda^+\cap\sigma}$ of non-zero scalars,
satisfying: $c_\lambda c_\mu=c_\nu$ for all
$\lambda,\mu,\nu\in\Lambda^+\cap\sigma$ such that 
$\End V_\lambda\cdot\End V_\mu$ contains $\End V_\nu$. Together with
the claim, this proves that any two structures of reductive semigroup
on $X$ are conjugate by a unique element of $\Aut^{\gxg}(X)$.
\end{proof}

Next we describe the idempotents in a self-adjoint stable
reductive semigroup $X$ (that is, the elements $e\in X$ such that
$e^2=e$), generalizing results of Putcha and Renner
\cite{Putcha82, Putcha84, Putcha_Renner88} about reductive
monoids. Clearly, the subset of idempotents is $\diag G$-invariant
and closed in $X$.

\begin{proposition}\label{idempotents}
Let $X$ be a self-adjoint stable reductive semigroup for $G$. Then
the set of idempotents of $X$ meets any $\gxg$-orbit $\cO$, along a
unique $\diag G$-orbit. Moreover, there exists a unique idempotent
$e_\cO$ in $\cO_0^{\diag T}$, where $\cO_0$ denotes the open
$B^-\times B$-orbit in $\cO$.
\end{proposition}

\begin{proof}
By Proposition \ref{irreducible} and Lemma \ref{square}, the closure
of $\cO$ is a self-adjoint reductive subsemigroup of $X$. Thus, we
may assume that $\cO$ is open in $X$.

Let $X'$ be the closure in $X$ of $\cO_0^{\diag T}$. By Lemma
\ref{square} and Theorem \ref{semigroup}, $X'$ is a subsemigroup of
$X$, which makes it a reductive semigroup (for $T$). Since $X'$ is an
affine toric variety, it follows easily that every $T\times T$-orbit
in $X'$ contains a unique idempotent. In particular, $\cO_0^{\diag T}$
contains a unique idempotent $e=e_\cO$.

Let $e'\in\cO$ be another idempotent. Then $e'=g e g'$ for
some $g,g'\in G$. Moving $e'$ in its $\diag G$-orbit, we may
assume that $e'= g e$ for some $g\in G$. Then $ege=e$, since
$e'{}^2=e'$. Let $H=\Stab_{\gxg}(e)$, then the structure of $H$ is given
by Proposition \ref{structure}. Since $g_1 e g_2^{-1} =e$ for any
$(g_1,g_2)\in H$, it follows that $g_1 e = e = e g_2$ for any 
$g_1\in R_u(P_I)H_J$ and $g_2\in R_u(P_I^-) H_J$; moreover, $xe=ex$
for any $x\in L_K$. By the Bruhat decomposition, one has 
$g = g_2 n x g_1$ for some $g_1\in R_u(P_I)$, $n\in N$, $x\in L_K$
and $g_2\in R_u(P_I^-) H_J$. Thus, 
$$
e = e g e = e n x e = e (nen^{-1}) n x.
$$
As a consequence, $e(nen^{-1})\in \cO^{\diag T}$. Since both
$e$, $nen^{-1}$ are idempotents in the stable toric semigroup 
$X^{\diag T}$, it follows that $nen^{-1}\in\cO^{\diag T}$. Thus, 
$n\in L_I$ by Proposition \ref{quasiaffine}; and since $L_I=H_J L_K$,
we may assume that $n=1$. Then $e = e x = x e$, and 
$$
e' = g e = g_2 x g_1 e = g_2 e = g_2 e g_2^{-1}.
$$
\end{proof}

For any idempotent $e$ in a semigroup $X$, the subset 
$$
eXe=\{x\in X~\vert~x=ex=xe\}
$$ 
is a submonoid, with identity element $e$. If, in addition, $X$ is a
self-adjoint stable reductive semigroup, then $eXe$ has a richer
structure: 

\begin{proposition}\label{submonoids}
Let $e$ be an idempotent in the self-adjoint stable reductive
semigroup $X$. Then the subsets 
$$
P=\{g\in G~\vert~ ge\in eG\},~ P^-=\{g\in G~\vert~ eg\in Ge\}
$$ are opposite parabolic subgroups of $G$, with common Levi subgroup 
$$
L=\{g\in G~\vert~ ge=eg\}\simeq \Stab_{\diag G}(e).
$$ 
Moreover, the subset $eXe$ is a normal reductive submonoid of $X$,
with identity element $e$ and unit group the quotient of 
$\Stab_{\diag G}(e)$ by a connected normal subgroup.
\end{proposition}

\begin{proof}
Note that 
$$
eXe\subseteq eXeG\subseteq eX\overline{GeG}\subseteq\overline{GeG},
$$
where the latter inclusion holds by Lemma \ref{square}. Thus,
$eXe=e\overline{GeG}e$, and we may assume that $GeG=\cO$ is open in $X$.
By Proposition \ref{idempotents}, we may also assume that 
$e=e_{\cO}$; then the isotropy group $\Stab_{\gxg}(e)$ is described by
Proposition \ref{structure}. Since $P$ (resp.~$P^-$) is the
projection of $\Stab_{\gxg}(e)$ to the first (resp.~second) copy of
$G$, this implies the first assertion.

Clearly, $eXe=\{x\in X~\vert~ x=xe=ex\}$ is closed in $X$, and
contains $eGeGe$ as a dense subset. Moreover, 
$R_u(P)LR_u(P^-)$ is dense in $G$, and $R_u(P)e=eR_u(P^-)=e$ by
Proposition \ref{structure}. Therefore, $eLeLe=eL=Le$ is dense in
$eXe$. In fact, the map $g\mapsto ge$ is a homomorphism of $L$ to the
unit group $G(eXe)$; hence the image $Le$ is open in that group. But
since $eXe$ is irreducible, $G(eXe)$ is connected, so that it
equals $Le$. In addition, the isotropy group 
$\Stab_L(e)=\{g\in G ~\vert~ ge=eg=e\}$ is connected by Proposition
\ref{structure} and Lemma \ref{automorphisms}.

It remains to show that the variety $eXe$ is normal. For this, we
argue as in the proof of the normality of $\overline{eTe}=X'$
(Proposition \ref{fixed}). We identify $k[eXe]$ with a 
$L\times L$-invariant subalgebra of $k[L/\Stab_L(e)]$; the latter
identifies in turn with a $L_K\times L_K$-invariant subalgebra of
$k[L_K]$, with the notation of Proposition \ref{structure}. Then
the vector space $k[eXe]$ is spanned by the functions on 
$L_K\simeq L_K\times L_K/\diag L_K$:
$$
(x,y)\mapsto\Tr_{V_\lambda}(uxp_{\lambda,K}y^{-1})
=\Tr_{V_\lambda}(p_{\lambda,K}y^{-1}ux),
$$
that is, by the matrix coefficients of $\End V_{\lambda,K}$
(where $\lambda\in\Lambda^+\cap\sigma$). It follows that 
$$
k[eXe] \simeq 
\bigoplus_{\lambda\in\Lambda^+\cap\sigma} \End V_{\lambda,K}
$$
as a $L_K\times L_K$-module. This implies normality of $eXe$, for
example by using Proposition \ref{connected} and Lemma \ref{normality}.
\end{proof}

Finally, we characterize self-adjoint stable reductive semigroups
among stable reductive varieties.

\begin{proposition}\label{stablesemired}
For a stable reductive variety $X=X_{\Sigma,t}$, the following
conditions are equivalent:

\begin{enumerate}

\item $X$ admits a structure of self-adjoint stable reductive
semigroup. 

\item The reference map $\rho$ is injective, and the cocycle $t$ is a
coboundary. 

\end{enumerate}

Then any two structures of stable reductive semigroup are conjugate by
a unique automorphism.
\end{proposition}

\begin{proof}
(1)$\Rightarrow$(2) Any $\sigma\in\Sigma$ defines a closed irreducible
$\gxg$-subvariety $Y_\sigma\subseteq X$. By Proposition
\ref{idempotents}, $Y_\sigma$ is equipped with a base point $e_\sigma$
in its open $B^-\times B$-orbit. This defines a unique isomorphism
$X_\sigma\simeq Y_\sigma$ of self-adjoint reductive semigroups, where
the semigroup structure on $X_\sigma$ is such that the base
point is idempotent (there is a unique such structure, by Theorem
\ref{semigroup} and Proposition \ref{idempotents}). These isomorphisms
identify the inclusions $Y_\tau\subseteq Y_\sigma$ with the maps
$i_{\tau\sigma}:X_\tau\to X_\sigma$. Thus, $X$ is isomorphic to the
direct limit of the $X_\sigma$, with the trivial cocycle
$t$. Moreover, $X$ is multiplicity-free by Lemma \ref{square}; it
follows that $\rho$ is injective.

(2)$\Rightarrow$(1) We may assume that $t$ is constant. Then the
preceding argument shows that $X=\varinjlim X_\sigma$ is a
self-adjoint semigroup satisfying conditions (1) and (2) of
Definition \ref{stablesemigroup}. Moreover, the multiplication
$m:X\times X\to X$ is surjective (since this holds on every
irreducible component). Therefore, the comorphism 
$$
m^{\#}:k[X]\to (k[X]\otimes_k k[X])^G
$$ 
is injective and $\gxg$-equivariant. On the other hand, the 
assumptions that $t$ is trivial and $\rho$ is constant imply that $X$
is multiplicity-free; it follows that $m^{\#}$ is an isomorphism.

The final assertion follows from Lemma \ref{square} and Theorem
\ref{semigroup}.
\end{proof}


\section{Families}

\subsection{General remarks} 

\label{subsection: Families} We record some general results on
arbitrary families of affine $G$-varieties, in the following sense.

\begin{definition}\label{families}
A \emph{family of affine (irreducible) $G$-varieties} over a scheme
$S$ is a scheme $\cX$ equipped with a morphism $\pi:\cX\to S$ and with
an action of the constant group scheme $G\times S$ over $S$, 
satisfying the following conditions:

\begin{enumerate}

\item $\pi$ is flat and affine. 

\item Every geometric fiber of $\pi$ is reduced (and irreducible).

\end{enumerate}

\end{definition}

This yields a rational action of $G\times S$ on the sheaf
$\pi_*\cO_{\cX}$, by  \cite{Mumford_GIT3ed} p.~25. As for families of
stable reductive varieties, we will often drop the adjective
``affine'' and deal with families of $G$-varieties.

\begin{lemma}\label{decomposition}
Let $\pi:\cX\to S$ be a family of multiplicity-finite 
(multi\-plicity-bounded, multi\-plicity-free) $G$-varieties.
Then the map $\cO_S\to (\pi_*\cO_{\cX})^G$ is an isomorphism, so that
$\pi$ is the categorical quotient by $G$. Moreover, one has an 
isomorphism of $\cO_S$-$G$-modules:
$$
\pi_*\cO_{\cX}\simeq\bigoplus_{\lambda\in\Lambda^+}
F_\lambda\otimes_k V_\lambda,
$$
where all the $F_\lambda\simeq (\pi_*\cO_{\cX})^U_\lambda$ are locally
free $\cO_S$-modules of finite rank (of bounded rank, of rank at
most $1$).

As a consequence, $\pi$ factors through a family $\pi//U:\cX//U\to S$ 
of multiplicity-finite (multiplicity-bounded,
multiplicity-free) $T$-varieties.
\end{lemma}

\begin{proof}
We may assume that $S$ is affine; hence so is $\cX$. Let
$\cR=\Gamma(\cX,\cO_{\cX})$ and $\cR_0=\Gamma(S,\cO_S)$. The
$\cR_0$-module $\cR$ is flat, so that its direct factor $\cR^G$ is
flat as well. Since every geometric fiber $\cX_{\os}$ is reduced and
multiplicity-finite, one has 
$$
\cR^G\otimes_{\cR_0} k(\os) \simeq
\Gamma(\cX_{\os},\cO_{\cX_{\os}})^G =
\Gamma(\cX_{\os},\cO_{\cX_{\os}})^{G(k(\os))} = k(\os).
$$
It follows
that the map $\cR_0 \to \cR^G$ is an isomorphism. And since $G$ acts 
rationally on $\cR$, one has an isomorphism of
$\cR_0$-$G$-modules: 
$$
\bigoplus_{\lambda\in\Lambda^+}
\Hom^G(V_\lambda,\cR)\otimes_k V_\lambda \to \cR
$$
given by evaluation. Since $\cR$ is a flat module over $\cR_0$,
the same holds for each $\Hom^G(V_\lambda,\cR)$. But this 
module of covariants is finitely generated over $\cR_0=\cR^G$ (since
$\cR$ is Noetherian), and hence locally free. Moreover, 
$\Hom^G(V_\lambda,\cR)\simeq \cR^U_\lambda$. The assertion on
ranks follows from the isomorphism 
$$
\cR^U_\lambda \otimes_{\cR_0} k(\os) \simeq 
\Gamma(\cX_{\os},\cO_{\cX_{\os}})^U_{\lambda} =
\Gamma(\cX_{\os},\cO_{\cX_{\os}})^{U(k(\os))}_{\lambda}.
$$
\end{proof}

Together with Proposition \ref{connected}, this implies at once the
following

\begin{corollary}\label{constant}
Let $\pi:\cX\to S$ be a family of $G$-varieties, where $S$ is
connected. Then the multiplicities of the $G(k(\os))$-module
$\Gamma(\cX_{\os},\cO_{\cX_{\os}})$ are independent of the
geometric point $\os$ of $S$. In particular, if some 
geometric fiber is multiplicity-finite (multiplicity-bounded,
multiplicity-free), then so are all geometric fibers. 

If, in addition, some geometric fiber is multiplicity-free with a
saturated weight set, then all geometric fibers are connected and
seminormal.
\end{corollary}

This defines the \emph{Hilbert function} 
$h:\Lambda^+\to\bN,~\lambda\mapsto \rank(F_{\lambda})$ 
of a family of multiplicity-finite $G$-varieties over a connected
scheme $S$.

Another useful observation is the following variant of a rigidity
lemma due to Knop (\cite{Knop94} Lemma 6.1.)

\begin{lemma}\label{loctriv}
Let $\pi:\cX\to S$ be a family of irreducible multiplicity-free
$G$-varieties. Then the family $\pi//U:\cX//U\to S$ is locally
trivial.
\end{lemma}

\begin{proof}
Choose $s\in S$. Replacing $S$ by a neighborhood of $s$, we may assume
that $S$, and hence $\cX$, is affine; replacing $\cX$ by $\cX//U$, we
may also assume that $G=T$ is a torus. Then
$$
\cR=\Gamma(\cX,\cO_{\cX})=\bigoplus_{\lambda\in\Lambda} \cR_\lambda,
$$ 
where each non-zero weight space $\cR_\lambda$ is a locally free 
module of rank $1$ over $\cR_0=\Gamma(S,\cO_S)$. Moreover,
$\cR_\lambda\neq 0$ if and only if $\lambda\in\Lambda^+_X$, where
$X=\cX_{\os}$. Replacing $T$ by a quotient, we may 
assume that its character group $\Lambda$ is generated by
$\Lambda_X^+$; then $T(k(\os))$ has an open orbit $\cO$ in $X$,
isomorphic to $T(k(\os))$. 

For any $\lambda$, $\mu$ in $\Lambda_X^+$, the multiplication map
$$
\cR_{\lambda}\otimes_{\cR_0} \cR_{\mu}\rightarrow \cR_{\lambda+\mu}
$$ 
induces an isomorphism on all geometric fibers; thus, this map is an
isomorphism. It follows that any local generator of the
$\cR_0$-module $\cR_{\lambda}$ is a non-zero divisor in $\cR$.

Choose $\lambda_0\in\Lambda^+_X$ in the interior of the cone $C_X$,
and let $f\in \cR_{\lambda_0}$ such that $f\vert_X\neq 0$. Then $f$
generates $\cR_{\lambda_0}$ at $s$, hence is a non-zero divisor in
$\cR$, and the open affine subset $\cX_f$ intersects $X$ along
$\cO$. As a consequence, 
$$
\cR[f^{-1}]=\bigoplus_{\lambda\in\Lambda} \cR[f^{-1}]_\lambda,
$$
where every $\cR[f^{-1}]_\lambda$ is a locally free $\cR_0$-module of
rank $1$. Choose a basis $(\lambda_1,\ldots,\lambda_r)$ of the lattice
$\Lambda$; for $1\leq i\leq r$, let $f_i\in \cR[f^{-1}]_{\lambda_i}$
such that $f_i\vert_X\neq 0$. Then, for any 
$\lambda=n_1\lambda_1+\cdots+n_r\lambda_r\in\Lambda$, the Laurent
monomial $f_{\lambda}=f_1^{n_1}\cdots f_r^{n_r}$ is a generator of the
$\cR_0$-module $\cR[f^{-1}]_\lambda$ at $s$. Thus, by shrinking $S$, we
may assume that $\cR[f^{-1}]$ is just the ring of Laurent polynomials
$\cR_0[f_1,f_1^{-1},\ldots,f_r,f_r^{-1}]
= \bigoplus_{\lambda\in\Lambda} \cR_0 f_\lambda$. 

Moreover, regarding $\cR$ as a $T$-invariant subalgebra of 
$\cR[f^{-1}]$, we have 
$$
\cR_{\lambda} = \cR[f^{-1}]_{\lambda}
$$
for all $\lambda\in \Lambda^+_X$. Indeed, $\cR[f^{-1}]_{\lambda}$ is
the union of the spaces $\cR_{\lambda+n\lambda_0}f^{-n}$ over all
non-negative integers $n$, and the equalities
$\cR_{\lambda+n\lambda_0} = \cR_{\lambda}\cR_{n\lambda_0}
= \cR_{\lambda} f^n$ yield
$\cR_{\lambda+n\lambda_0}f^{-n}=\cR_\lambda$. 

So we have constructed generators $f_\lambda$ of the
$\cR_\lambda$'s, satisfying $f_\lambda f_\mu = f_{\lambda+\mu}$. This
trivializes $\pi$ in a neighborhood of $s$.
\end{proof}

Finally, we show how to construct families of $G$-varieties with a
prescribed fiber, starting from certain families of 
$\gxg$-varieties; this will be developed in Subsection 
\ref{subsec: The Vinberg family}. 

\begin{lemma}\label{construction}
Let $X$ be a $G$-variety and let $\pi:\cX\to S$ be a family of
$\gxg$-varieties over a connected scheme, such that 
$\cX_{s_0}\simeq G$ for some $s_0\in S$ with $k(s_0)=k$. Then the
categorical quotient
$$
X * \cX=(X\times \cX)//G
$$ 
(for the action by $g\cdot(x,\xi)=(gx,(g,1)\xi)$), equipped with the
$G$-action induced by the $\{1\}\times G$-action on $\cX$, and with
the map
$$
\pi_X:(X\times \cX)//G \to \cX//G \to S,
$$ 
is a family of $G$-varieties, with fiber at $s_0$ isomorphic to $X$.

If, in addition, $X$ is multiplicity-finite
(multiplicity-bounded, multiplicity-free), then so are all
geometric fibers of $\pi_X$.
\end{lemma}

\begin{proof}
We may assume again that $S$, and hence $\cX$, is affine. Let again
$\cR=\Gamma(\cX,\cO_{\cX})$ and $\cR_0=\Gamma(S,\cO_S)$. Then 
$$
\cR=\bigoplus_{\lambda\in\Lambda^+} 
F_\lambda\otimes_k V_\lambda^*\otimes_k V_\lambda
$$
as a $\cR_0$-$\gxg$-module, where every $F_\lambda$ is an
invertible $\cR_0$-module (this follows from Lemma
\ref{decomposition} together with the decomposition  
$$
k[G]\simeq\bigoplus_{\lambda\in\Lambda^+} V_\lambda^*\otimes_k V_\lambda
$$ 
as a $\gxg$-module.) We also have 
$$
k[X]=\bigoplus_{\lambda\in\Lambda^+} 
\Hom^G(V_\lambda,k[X])\otimes_k V_\lambda
$$
as a $G$-module. It follows that
$$
\Gamma(X * \cX,\cO_{X * \cX}) = (k[X]\otimes_k R)^G
= \bigoplus_{\lambda\in\Lambda^+}
\Hom^G(V_\lambda,k[X])\otimes_k F_\lambda\otimes_k V_\lambda
$$
as a $\cR_0$-$G$-module, where $\cR_0$ acts on $F_\lambda$, and
$G$ acts on $V_\lambda$. Thus, $\pi_X$ is flat. Moreover, 
$(X * \cX)_{\os} \simeq (X\times \cX_{\os})//G(k(\os))$ for any 
geometric point $\os$ of $S$, so that the geometric fibers of $\pi_X$
are reduced. Finally, if $\cX_{s_0}$ is isomorphic to
$G$, then $(X * \cX)_{s_0}\simeq (X\times G)//G \simeq X$. This proves
the first assertion.

If, in addition, $X$ is multiplicity-finite
(multiplicity-bounded, multiplicity-free), then the dimensions
of all spaces $\Hom^G(V_\lambda,k[X])$ are finite (bounded, at
most $1$). This implies the second assertion.
\end{proof}

\subsection{Moduli of embedded stable reductive varieties}

Next we show that the families of multiplicity-finite
$G$-subvarieties of a fixed $G$-module $V$, with
a fixed Hilbert function $h$, admit a fine moduli space. In the case
where $G$ is abelian, $G$-modules are nothing but $k$-vector spaces
endowed with a grading by the character group of $G$. In that case,
existence of the moduli space follows from the construction by Haiman
and Sturmfels \cite{HaimanSturmfels} of the multigraded Hilbert
scheme, parametrizing all homogeneous quotient rings of a graded
polynomial ring with a fixed Hilbert function.

As a generalization of this setting to arbitrary $G$, consider a
function $h:\Lambda^+\to\bN$ and a finite-dimensional $G$-module
$V$. Define a functor
$$
\cM_{h,V}:\operatorname{(Schemes)}^o \to \operatorname{(Sets)}
$$
by assigning to $S$ the set of all multiplicity-finite subfamilies
$\pi:\cX\to S$ of $(V\times S\to S)$, with Hilbert function $h$.

\begin{theorem}\label{Hilb}
The functor $\cM_{h,V}$ is representable by a scheme $M_{h,V}$ which
is quasiprojective over $k$.
\end{theorem}

\begin{proof}
Note that we may define families of affine $G$-schemes by dropping
condition (2) (on reducedness of geometric fibers) in Definition
\ref{families}; the notions of multiplicity-finiteness and Hilbert
function adapt easily to such families. Thus, we may consider the
functor  
$$
\cH_{h,V}:\operatorname{(Schemes)}^o \to \operatorname{(Sets)}
$$
by assigning to $S$ the set of all closed $G$-invariant subschemes
$\cX\subseteq V\times S$, flat and multiplicity-finite over $S$ with
Hilbert function $h$.
We will prove that $\cH_{h,V}$ is representable by a quasi-projective
scheme $H_{h,V}$, and then that $\cM_{h,V}$ is representable by an open
subscheme of $H_{h,V}$.

For this, we reduce to a situation where the main result of Haiman and
Sturmfels applies, as follows. Choose a regular dominant
one-parameter subgroup $\gamma:\bG_m\to T$. Then the evaluation map
$$
\Lambda\to \bZ,~\lambda\mapsto \langle\lambda,\gamma\rangle
$$ 
satisfies the following properties:

\begin{enumerate}

\item
$\langle\lambda,\gamma\rangle\ge 0$ for all $\lambda\in\Lambda^+$.

\item
$\langle\nu,\gamma\rangle \le 
\langle\lambda,\gamma\rangle + \langle\mu,\gamma\rangle$
for all $\lambda$, $\mu$, $\nu$ in $\Lambda^+$ such that the
$G$-module $V_{\nu}$ occurs in the decomposition of 
$V_{\lambda}\otimes_k V_{\mu}$.

\item
For any non-negative integer $n$ and any character $\chi$ of the
center $Z(G)$, the set 
$$
\{\lambda\in\Lambda^+~\vert~ \langle\lambda,\gamma\rangle = n
\text{ and } \lambda\vert_{Z(G)}=\chi\}
$$
is finite.

\end{enumerate}

Now consider an affine scheme $X$ with a $G$-action, and the
corresponding $G$-algebra 
$$
R = k[X] = \bigoplus_{\lambda\in\Lambda^+} 
\Hom^G(V_{\lambda}, R)\otimes_k V_{\lambda}.$$ 
For any non-negative integer $n$, let
$$
R_{\le n}= \bigoplus_{\lambda\in\Lambda^+,
\langle\lambda,\gamma\rangle \le n} 
\Hom^G(V_{\lambda}, R)\otimes_k V_{\lambda}.
$$
This yields an ascending filtration of $R$ by $G$-submodules, and (2)
shows that it is an algebra filtration. Let $z$ be a variable, and let
$$
\cR = \bigoplus_{n\ge 0} R_{\le n} z^n
$$
be the corresponding Rees algebra. This is a graded $G$-invariant
subalgebra of $R[z]$, containing $k[z]$. Any $G$-invariant ideal
$I\subseteq R$ yields a graded $G$-invariant ideal
$$
\cI = \bigoplus_{n\ge 0} (I\cap R_{\le n}) z^n \subseteq \cR,
$$
and hence a graded $G$-algebra $\bar{\cR} = \cR/\cI$.
One checks easily that the assignment $I\mapsto\cI$ sets up a
bijective correspondence from $G$-invariant ideals of $R$, to those
graded $G$-invariant ideals $\cI$ of $\cR$ such that $z$ is a
non-zero divisor in $\cR/\cI$. Moreover, $I$ is prime (radical) if
and only if so is $\cI$. Finally, by (3), the $G$-algebra $R/I$ is
multiplicity-finite if and only if so is the 
$\bG_m\times Z(G)$-algebra $\cR/\cI$, where the action of $\bG_m$ is
defined by the compatible gradings of $\cR$ and $\cI$. 

Another graded $G$-invariant ideal of $\cR$ is 
$$
\cJ = \cI+z\cR = \bigoplus_{n=0}^{\infty} 
(R_{\le n-1} + I\cap R_{\le n}) z^n.
$$
Note the isomorphisms of $G$-modules 
$$
\cR_n/\cJ_n \simeq R_{=n}/I_{=n} \simeq 
\bar{\cR}_{\le n}/\bar{\cR}_{\le n-1}
$$ 
with evident notation. 
Thus, the Hilbert functions of $R/I$ (as a $G$-algebra),
of $\cR/\cI$ and of $\cR/\cJ$ (as $\bG_m\times G$-algebras) determine
each other. In particular, the Hilbert function of $R/I$ as a
$G$-algebra determines that of $\cR/\cI$ and $\cR/\cJ$ as
$\bG_m\times Z(G)$-algebras.

Next let $\cX\subseteq V\times S$ be a $G$-invariant closed
subscheme. Let 
$$
R=\Sym^*(V^*)\otimes_k \cO_S
$$ 
be the direct image of the structure sheaf
$\cO_{V\times S}$ under projection to $S$. Then $\cX$ yields a
$G$-invariant sheaf of ideals $I\subseteq R$, and hence a
$\bG_m\times G$-invariant sheaf of ideals $\cI$ in the
Rees algebra
$$
\cR=(\bigoplus_{n=0}^{\infty} \Sym^*(V^*)_{\le n} z^n)
\otimes_k \cO_S.
$$
The argument of Lemma \ref{decomposition} shows that
$\cX$ is flat over $S$ if and only if so is $\cR/\cI$. Then
$\cJ=\cI+z\cR$ is another $\bG_m\times G$-invariant sheaf of ideals,
and $\cR/\cJ$ is flat over $S$, too.

By \cite{HaimanSturmfels} Theorem 1.1, the families of pairs 
$(\cR/\cI,\cR/\cJ)$ of $\bG_m\times Z(G)$-invariant quotients
of $\cR$ with fixed Hilbert functions (determined by $h$) are
parametrized by a quasi-projective scheme $S_{h,V}$. The
$\cO_{S_{h,V}}$-modules $\cR/\cI$ and $\cR/\cJ$ split up into direct
sums of locally free modules of finite rank according to the action of 
$\bG_m\times Z(G)$. By Lemmas~\ref{lem:univ1} and \ref{lem:univ2}
below, those pairs where $\cI+z\cR = \cJ$ and both $\cI$, $\cJ$ are
$G$-invariant are parametrized by a closed subscheme $H_{h,V}$ of
$S_{h,V}$.

We claim that $H_{h,V}$ represents the functor $\cH_{h,V}$. By the
preceding discussion, it suffices to check that $z$ is a non-zero
divisor in $\cR/\cI$. Equivalently, the multiplication by $z$ fits
into exact sequences
$$
0 \to \cR_{n-1}/\cI_{n-1} \to \cR_n/\cI_n \to \cR_n/\cJ_n \to 0
$$
for all $n$. But these sequences are right exact and split up (under
the $\bG_m$-action) into direct sums of sequences involving only
locally free $\cO_{H_{h,V}}$-modules of finite rank. Hence their
exactness is equivalent to a collection of numerical conditions on
these ranks, i.e., it only depends on the Hilbert function. This
proves our claim. 

It remains to show that $\cM_{h,V}$ is representable by an open
subscheme of $H_{h,V}$. For this, we consider the universal family
$\pi:U_{h,V}\to H_{h,V}$. We are interested in the locus of $H_{h,V}$
where the fibers of $\pi$ are reduced, and we know that the preimage
of this locus in $U_{h,V}$ is open and invariant under the natural
$G$-action. To conclude, note that $\pi$ maps any closed
$G$-invariant subscheme to a closed subscheme, since it factors as
the quotient $U_{h,V}\to U_{h,V}//G$ followed by a finite 
morphism $U_{h,V}//G \to H_{h,V}$. 
\end{proof}

\begin{lemma}\label{lem:univ1}
  Let $S$ be a scheme, $V$ a $k$-vector space, and let $F_1,F_2$ be
  two subsheaves of $V\otimes_k \cO_S$ such that the quotients
  $V\otimes_k \cO_S/F_1=Q_1$, $V\otimes_k \cO_S/F_2=Q_2$ are 
  locally free of finite rank. Then there exists a closed subscheme
  $S(F_1,F_2)\subseteq S$ satisfying the following universal
  property: for any morphism $f:S'\to S$, $f^*F_1$ is a subsheaf of
  $f^*F_2$ if and only if $f$ factors through $S(F_1,F_2)$.
\end{lemma}
\begin{proof}
  By shrinking, we can assume $S$ to be affine, $S=\Spec R$, and
  $Q_1,Q_2$ to be freely generated by the images of certain vectors of
  $V$. Then each $F_i$ contains $V_i\otimes_k \cO_S$ for some
  subspace $V_i\subseteq V$ of finite codimension. Replacing $V$ with 
  $V/V_1\cap V_2$, we can assume that $V$ is finite dimensional.
  Shrinking again, we can also assume that both $F_1$, $F_2$ are free
  with bases $\{a_1,\dots, a_{n_1}\}$, $\{b_1,\dots, b_{n_2}\}$
  respectively. The second basis can be extended to a basis 
  $\{b_1,\dots, b_n\}$ of $V\otimes_k \cO_S$. Write the basis of $F_1$
  as an $n\times n_1$ matrix $(A_{ij})$ with entries in $R$. Then
  $S(F_1,F_2)$ is the zero scheme of entries $A_{ij}$ with $i> n_2$. 
\end{proof}

\begin{lemma}\label{lem:univ2}
  Let $S$ be a scheme, $V$ a rational $G$-module, and 
  $F\subseteq V\otimes_k\cO_S$ a subsheaf such that the quotient
  $V\otimes_k\cO_S/F$ is locally free of finite rank. Then there
  exists a closed subscheme $S_F\subseteq S$ satisfying  
  the following universal property: for any morphism $f:S'\to S$,
  $f^*F$ is a $G$-invariant subsheaf of 
  $V\otimes_k\cO_{S'}$ if and only if $f$ factors through $S_F$.
\end{lemma}
\begin{proof}
  Indeed, $S_F$ is the intersection of the closed subschemes
  $S(F,g^*F)$ for $g\in G$, which exist by the previous Lemma.
\end{proof}


\subsection{Local isotriviality}

An important consequence of the classification of reductive varieties
(Theorem \ref{correspondence1}) is the following 

\begin{corollary}\label{isotrivial}
Let $\pi:\cX\to S$ be a family of reductive varieties over an integral
scheme. Then there exists an \'etale morphism $\varphi:S'\to S$ such
that the pull-back family $\cX\times_S S'\to S'$ is trivial. In
particular, $\pi$ admits a general fiber.
\end{corollary}

\begin{proof}
We may assume that $S$, and hence $\cX$, are affine. 
Let $\cR=\Gamma(\cX,\cO_{\cX})$ and $\cR_0=\Gamma(S,\cO_S)$. Write 
$$
\cR=\bigoplus_{\lambda\in\Lambda^+} F_\lambda\otimes_k \End V_\lambda
$$
as a $\cR_0$-$\gxg$-module; then $\cR_0 = F_0$. By shrinking $S$, we
may further assume that any non-zero $F_\lambda$ is generated by
some $f_{\lambda}$, satisfying 
$f_{\lambda} \cdot  f_{\mu} = f_{\lambda+\mu}$ (Lemma \ref{loctriv}).

On the other hand, the geometric generic fiber of $\pi$ is a reductive
variety over $k(\oS)$; let $C$, $K$ be the corresponding combinatorial
data. Then
$$
k(\oS)\otimes_{\cR_0}\cR = \bigoplus_{\lambda\in\Lambda\cap C}
k(\oS)\otimes_k \End V_\lambda = k(\oS)\otimes_k k[X_{C,K}]
$$
as a $k(\oS)$-algebra. It follows that each $F_\lambda$ is contained
in $k(\oS)$, and is nonzero if and only if $\lambda\in\Lambda\cap C$.
This yields elements $f_\lambda\in k(\oS)^*$, which are contained in
some finite extension $L$ of $k(S)$ by the first step of the
proof. So we obtain an algebra isomorphism 
$$
L\otimes_{\cR_0} \cR =  L\otimes_k k[X_{C,K}]
$$
which implies our statement.
\end{proof}

\subsection{One-parameter degenerations}

\label{subsec: One-parameter degenerations}

We construct certain families of stable reductive varieties over the
affine line; these will provide local models for all one-parameter
families of stable reductive varieties. We begin with families
whose generic fiber is geometrically irreducible. 

\begin{definition}\label{function}
A \emph{$W$-admissible height function} is a map $h:\sigma\to\bR$, where
$\sigma$ is a $W$-admissible cone in $\Lambda_{\bR}$, satisfying the
following conditions:

\begin{enumerate}

\item (on rationality) $h(\lambda)\in\bQ$ for any 
$\lambda\in\sigma\cap\Lambda$. 

\item (on invariance) $h(wx)=h(x)$ for any $w\in\Stab_W(\sigma)$ and
$x\in\sigma$. 

\item (on convexity) $h(x_1+x_2)\leq h(x_1)+h(x_2)$ for all $x_1$,
$x_2$ in $\sigma$.

\item (on piecewise linearity) $h$ is linear on each cone of some
finite subdivision of $\sigma$.

\end{enumerate}

In particular, $h$ is positively homogeneous.
\end{definition}

We then denote by $\Sigma_h$ the set of cones in $\Lambda_{\bR}$
consisting of all $W$-translates of the maximal cones where $h$ is
linear, together with their faces. Then the union $\vert\Sigma_h\vert$
of these cones is just $\bigcup_{w\in W} w\sigma$.
Clearly, $\Sigma_h$ is a $W$-complex of cones, the reference map 
$\vert\Sigma_h\vert\to \Lambda_{\bR}$ being the inclusion.

\begin{definition}\label{limit}
Let $X_h$ be the stable reductive variety associated with the
$W$-complex of cones $\Sigma_h$ and with the trivial cocycle $t$.
\end{definition}

Note that $X_h$ is unchanged when replacing $h$ by a positive rational
multiple.

We now construct a one-parameter family of stable reductive varieties
with general fiber $X_\sigma$ and special fiber $X_h$.
Let $\tilde G=G\times \bG_m$; this is a connected reductive group with
maximal torus $\tilde T=T\times \bG_m$, weight lattice
$\tilde\Lambda=\Lambda\times\bZ$, and Weyl group isomorphic to $W$. Let 
$$
\tilde \sigma_h=\{(x,y)\in\tilde\Lambda_{\bR}=\Lambda_{\bR}\times\bR
~\vert~ x\in\sigma,~ y\in \bR,~ h(x)\leq y\}.
$$
One easily checks that $\tilde\sigma_h$ is a $W$-admissible cone in
$\tilde\Lambda_{\bR}$; let $\tilde X_h$ be the corresponding reductive
variety (for $\tilde G$). Since $\tilde \sigma_h$ contains 
$\{0\}\times\bR_{\geq 0}$, we obtain a non-zero map
$$
\pi_h:\tilde X_h \to \bA^1
$$
which is $\gxg$-invariant and $\bG_m\times \bG_m$-equivariant 
(where $\bG_m \times \bG_m$ acts on $\bA^1$ by
$(t_1,t_2)z = t_1t_2^{-1}z$).

\begin{proposition}\label{degeneration}
$\pi_h$ is surjective, and 
$\pi_h^{-1}(\bA^1 - \{0\})\simeq X_\sigma\times(\bA^1 - \{0\})$ 
as $\tilde G\times\tilde G$-varieties.

Moreover, the fiber of $\pi_h$ at $0$ is reduced if and only if: 
$h(\lambda)\in\bZ$ for any $\lambda\in\Lambda\cap\sigma$. Then
$\pi_h$ is a flat family of stable reductive varieties, and 
its fiber at $0$ is isomorphic to $X_h$.
\end{proposition}

\begin{proof}
Let $K=K_\sigma$ and $\cX=\tilde X_h$. Write 
$$
k[\tilde G]_{(K)}=k[G]_{(K)}[z,z^{-1}]\simeq
\bigoplus_{\lambda\in\Lambda^+,\,n\in\bZ} z^n\End V_\lambda,
$$ 
where $z$ is the coordinate function on $\bG_m$. Then 
$$
k[\cX]\simeq
\bigoplus_{\lambda\in\Lambda^+\cap\sigma,\,n\in\bZ,\,h(\lambda)\leq n}
z^n\End V_\lambda
$$
as a subalgebra of $k[G]_{(K)}[z,z^{-1}]$, and the map 
$\pi_h:\cX \to \bA^1$ corresponds to the inclusion  
$k[z]\subseteq k[\cX]$. Moreover, one has
$k[\cX][z^{-1}]=k[X_\sigma][z,z^{-1}]$. This trivializes
$\pi_h$ over $\bA^1 - \{0\}$.

For any $n\in\bZ$, let
$$
k[X_\sigma]_{\leq n}=
\bigoplus_{\lambda\in\Lambda^+\cap\sigma,\,h(\lambda)\leq n}
\End V_\lambda.
$$
Then the $k[X_\sigma]_{\leq n}$ define an ascending filtration of
the algebra $k[X_\sigma]$ (indeed, if $\lambda$, $\mu$, $\nu$ are in
$\Lambda^+\cap\sigma$ and satisfy $\nu\leq_K\lambda+\mu$, then 
$\nu\in\Conv(W_K\lambda)+\Conv(W_K\mu)$, whence 
$h(\nu)\leq h(\lambda)+h(\mu)$ by invariance and convexity of $h$).
Moreover, $k[\cX]$ is the Rees algebra associated with this
filtration. It follows that
$$
k[\cX_0] = k[\cX]/zk[\cX] \simeq
\bigoplus_{\lambda\in\Lambda^+\cap\sigma} \End V_\lambda,
$$
as a $\gxg$-module, the multiplication being given by
$$
f_{\lambda}\cdot f_{\mu} = 
\sum_{\nu\in\Lambda^+,\,\nu\leq_K\lambda+\mu,\,
\ulcorner h(\nu)\urcorner  = \ulcorner h(\lambda) \urcorner
+ \ulcorner h(\mu) \urcorner} 
\Tr_1(p^{\nu}_{\lambda\mu}(f_{\lambda}\otimes_k f_{\mu}))
$$
for any $f_\lambda\in \End V_\lambda$ and $f_\mu\in\End V_\mu$. Here
for any $y\in \bR$, we denote by $\ulcorner y\urcorner$ the largest
integer $n$ such that $n\leq y$.

As a consequence, if there exists $\lambda\in \Lambda^+\cap\sigma$
such that $h(\lambda)\notin\bZ$, then any $U^-\times U$-invariant
function $f_\lambda\in\End V_\lambda$ is nilpotent; thus,
$\cX_0$ is not reduced. On the other hand, if
$h(\Lambda^+\cap\sigma)\subseteq \bZ$, then one has for any 
$U^-\times U$-invariant functions $f_\lambda$, $f_\mu$: 
$$
f_{\lambda}\cdot f_{\mu}\neq 0 \text{ if and only if }
h(\lambda)+h(\mu) = h(\lambda+\mu).
$$
Since $h$ is positively homogeneous, it follows that the algebra
$k[\cX_0]^{U^-\times U}$ is reduced, so that $\cX_0$
is reduced by Lemma \ref{Uinv}. Moreover, its irreducible components are
the $X_\tau$, where $\tau\subseteq\sigma$ is a maximal cone such that:
$h\vert_{\tau}$ is linear, and $\tau^0$ meets $\Lambda^+_{\bR}$. Since
$\cX_0$ is a stable reductive variety (by Corollary \ref{constant}),
it follows easily that it is isomorphic to $X_h$. 
\end{proof}

Next we consider a family $\pi:\cX\to S$ of stable reductive varieties
over a nonsingular curve, that is, $S$ is integral, regular and of
dimension $1$. We assume that the generic fiber of $\pi$ is
geometrically irreducible; let $\sigma$ be the corresponding
$W$-admissible cone. We show that $\pi$ looks locally like a family
associated to a $W$-admissible height function on $\sigma$. 

\begin{proposition}\label{curve}
With the preceding notation, for any closed point $s\in S$ there exist
a nonsingular curve $S'$, a closed point $s'\in S'$,
two non-constant morphisms $\varphi:(S',s')\to (S,s)$ and
$z:(S',s')\to(\bA^1,0)$, and a unique $W$-admissible height function
$h:\sigma\to\bR$, such that the pull-back families $\cX\times_S S'$
and $\tilde X_h\times_{\bA^1} S'$ are isomorphic. In particular, the
fiber $\cX_s$ is isomorphic to $X_h$.
\end{proposition}

\begin{proof}
Using Corollary \ref{isotrivial}, we may assume that $\pi$ is trivial
over the punctured curve $S-\{s\}$. We may further assume that $S$ is
affine, regular, and admits a morphism $z : S \to \bA^1$ having a
unique zero at $s$, of order $1$. Then $\cX$ is affine. Setting
$\cR = \Gamma(\cX,\cO_{\cX})$ and $\cR_0=\Gamma(S,\cO_S)$, we have 
$\cR_0[1/z]=\Gamma(S-\{s\},\cO_S)$. Moreover, Lemma \ref{decomposition}
yields an isomorphism of $\cR_0$-$\gxg$-modules
$$
\cR =  \bigoplus_{\lambda\in\Lambda^+\cap\sigma} 
F_\lambda\otimes_k \End V_\lambda,
$$
embedded in
$$
\Gamma(\cX-\cX_s,\cO_{\cX}) = \cR_0[1/z]\otimes_k k[X_\sigma]
=\bigoplus_{\lambda\in\Lambda^+\cap\sigma} 
\cR_0[1/z]\otimes_k \End V_\lambda
$$
as a subalgebra. Here each $F_\lambda$ is an invertible
$\cR_0$-submodule of $\cR_0[1/z]$, trivial over $\cR_0[1/z]$. It
follows that $F_\lambda=z^{h(\lambda)} \cR_0$ for a unique map
$h:\Lambda^+\cap\sigma\to \bZ$. We have for any $\lambda$, $\mu$ in
$\Lambda^+\cap\sigma$:
$$
z^{h(\lambda)}\End V_\lambda \cdot z^{h(\mu)} \End V_\mu =
\bigoplus_{\nu\in\Lambda^+,\,\nu\leq_K \lambda+\mu}
z^{h(\lambda)+h(\mu)} \End V_\nu,
$$
where $K=K_\sigma$. It follows that 
$h(\nu)\leq h(\lambda) + h(\mu)$ for all such triples $\lambda$,
$\mu$, $\nu$, and hence that the subspace
$$
R=\bigoplus_{\lambda\in\Lambda^+\cap C,\,n\in\bZ,\,h(\lambda)\leq n}
z^n \End V_\lambda\subseteq k[X_\sigma][z,z^{-1}]\subseteq \cR
$$
is in fact a $\gxg$-invariant subalgebra containing $k[z]$.

Clearly, the natural map $R\otimes_{k[z]}\cR_0\to \cR$ is an
algebra isomorphism. As a consequence, the map 
$R/zR\to k[\cX_s]$ is an isomorphism; in particular,
the algebra $R/zR$ is finitely generated, so that the same holds for 
$R^{U^-\times U}/zR^{U^-\times U}$. Now one easily checks that the
algebra $R^{U^-\times U}$ is generated by $z$, together with any family
of $T\times T$-eigenvectors which maps to generators of the
quotient $R^{U^-\times U}/zR^{U^-\times U}$. It follows that $R$ is
finitely generated.

Let $\tilde X=\Spec R$. Then $\tilde X$ is equipped with a
$\gxg$-action, and also with a compatible $\bG_m$-action (since $R$
is graded by the degree in $z$). Thus, $\tilde X$ is a 
$\tilde G\times\tilde G$-variety, which is clearly irreducible,
diagonal, and multiplicity-free, with weight set 
$$
\{(\lambda,n) ~\vert~ \lambda\in\Lambda^+\cap\sigma,
~n\in\bZ,~h(\lambda)\leq n\}.
$$
We claim that this set is saturated in
$\tilde\Lambda = \Lambda\times\bZ$. Indeed, since $R/zR$ is 
reduced, one obtains $h(m\lambda) = m h(\lambda)$ for all 
$m\in \bZ_{\geq 0}$ and $\lambda\in\Lambda^+\cap\sigma$, by arguing as
in the proof of Proposition \ref{degeneration}. This implies the claim.

By that claim and Proposition \ref{irreducible}, 
$\tilde X$ is a reductive variety for $\tilde G$. Now one easily
checks that $\tilde X=\tilde X_h$ for a unique $W$-admissible height
function $h:\sigma\to\bR$ extending our original
$h:\Lambda^+\cap\sigma\to\bZ$.
\end{proof}

We can now tackle the general case, when the generic fiber of a family 
is not geometrically irreducible. Let $\Sigma$ be a $W$-complex of
cones.

\begin{definition}\label{system-of-functions}
A \emph{$W$-admissible system of heights} on $\Sigma$ is a
collection of height functions 
$\uh= \{h_{\sigma}:\sigma\to\bR \,|\,\sigma\in \Sigma \}$ 
such that for any $\sigma_1,\sigma_2,\tau$ in $\Sigma$ with 
$\tau \prec \sigma_1\cap \sigma_2$, the difference
$(h_{\sigma_1} - h_{\sigma_2})|_{\tau}$ is a linear function.
\end{definition}

Not let $t\in Z^1(\Sigma, \Aut)$ be a 1-cocycle defined over the
field $k(z)$. We can write $t=z^\gamma t'$, where $t'$ is a 1-cocycle
which is invertible at $z=0$, and $\gamma$ is obtained by taking the
valuation of $t$ at $0$. Explicitly, $\gamma$ is a collection of 
homomorphisms 
$\gamma_{\tau}\in \Hom(\Lambda\cap\lin{\tau}/K_{\tau}, \bZ)$ 
for all $\tau\prec \sigma_1\cap\sigma_2$
(that is, of one-parameter subgroups of $\Aut^{\gxg}(X_\tau)$)
which satisfy the 1-cocycle condition on all
triple intersections. We call $t$ admissible if the following
condition holds for all faces $\tau\prec\sigma_1\cap \sigma_2$:
$$
(h_{\sigma_1} - h_{\sigma_2})|_{\tau} = \gamma_{\tau}.
$$
Then we can define a family $\tilde X_{\uh,t}$ 
over the open subset of $\bA^1$ (containing $0$) where $t'$ is
regular. It will be a family of stable reductive varieties precisely
when all $h_{\sigma}$ are $\bZ$-valued.

Indeed, each height function $h_{\sigma}$ defines a
subdivision of the cone $\sigma$ and a family $\tilde X_{h_{\sigma}}$
corresponding to a subalgebra of $k[G][z,z^{-1}]$. We then just glue
them along $t$. It is clear that the general fiber of
$\tilde X_{\uh, t}$ is the union of varieties
$X_{h_{\sigma}}$ glued along~$t'(0)$. Vice versa, one has

\begin{proposition}\label{prop:curve-stable}
  Let $\pi:\cX \to S$ be a family of stable reductive varieties over a
  nonsingular curve germ $(S,s)$, and assume that the generic fiber
  is multiplicity-free. Then after a nonconstant base change
  $(S',s')\to (S,s)$ the family $\cX'$ is a pullback of some 
  $\tilde X_{\uh, t}$ under $(S',s')\to (\bA^1,0)$.
\end{proposition}

\begin{proof}
  After making a finite base change, we can assume that every
  irreducible component $\cY_{\eta}$ of the generic fiber
  $\cX_{\eta}$ is geometrically irreducible and is trivial over
  $k(\eta)$. Its closure $\overline{\cY_\eta}=\cY$ in $\cX$ is
  $\gxg$-invariant and reduced, so must be flat over $S$.  Thus, $\cY$
  is a family of stable reductive varieties, with geometrically
  irreducible generic fiber. By Proposition \ref{curve}, we may assume
  that $\cY$ is a pull-back of $\tilde X_{h_{\sigma}}$, where
  $h_\sigma$ is a well-defined weight function. This defines a
  $W$-complex of cones with an injective reference map, since $\cX$
  is multiplicity-free. Any intersection $\tau=\sigma_1\cap\sigma_2$
  is a cone, and $\tilde X_{h_{\tau}}$ embeds into both
  $\tilde X_{h_{\sigma_1}}$ and $\tilde X_{h_{\sigma_2}}$. These
  identifications define an admissible cocycle $t$. 
\end{proof}


\subsection{The Vinberg family}
\label{subsec: The Vinberg family}

We construct a family of $\gxg$-varieties with general fiber $G$,
after Vinberg \cite{Vinberg95}. Consider the algebra 
$$
k[T\times G]=\bigoplus_{\lambda\in\Lambda^+,\,\mu\in\Lambda}
e^{\mu}\End V_\lambda
$$
(where $e^{\mu}$ denotes the weight $\mu$ regarded in $k[T]$),
and its subspace
$$
R=\bigoplus_{\lambda\in\Lambda^+,\,\mu\in\Lambda,\,\lambda\leq_{\Pi}\mu}
e^{\mu}\End V_\lambda.
$$
Then $R$ is invariant under the product action of $T\times\gxg$,
and is also a subring of $k[T\times G]$, since 
$$
\End V_{\lambda}\cdot\End V_{\mu}\subseteq
\bigoplus_{\nu\leq_{\Pi}\lambda+\mu}\End V_\nu
$$
for all $\lambda,\mu\in\Lambda^+$. 
Note that the invariant subalgebra
$$
R_0 = R^{\gxg} \simeq 
\bigoplus_{\mu\in\Lambda,\,\mu\geq_{\Pi}0} k e^{\mu}
$$
is the polynomial ring $k[e^{\alpha}, \alpha\in\Pi]$, that is, the
coordinate ring of the affine space $\bA^{\Pi}$.
Moreover, the subalgebra
$$
R^{U^-\times U}=\bigoplus_{\lambda\in\Lambda^+,\,\mu\in\Lambda}
e^{\mu}(\End V_\lambda)^{U^-\times U}
$$
is isomorphic to $k[\Lambda^+] \otimes_k R_0$,
where $k[\Lambda^+]=k[G]^{U^-\times U}$ is the algebra of the
monoid of dominant weights. In particular, the algebra 
$R^{U^-\times U}$ is a finitely generated normal domain; thus, the
same holds for $R$ by \cite{Grosshans97} Theorem 16.2. So
$\cV=\Spec R$ is a $T\times\gxg$-variety, endowed with a morphism
$$
p:\cV\to \bA^{\Pi}
$$
that is $\gxg$-invariant and $T$-equivariant (where $T$ acts
linearly on $\bA^{\Pi}$ with weights the simple roots). One has 
an isomorphism of $R_0$-modules
$$
R \simeq \bigoplus_{\lambda\in\Lambda^+} 
R_0 \otimes_k e^{\lambda}\End V_\lambda.
$$
It follows that $p$ is flat, and that 
$p//(U^-\times U):\cV//(U^-\times U)\to \bA^{\Pi}$ is a trivial family
of $T$-varieties; its fiber is the affine toric variety associated
with the monoid $\Lambda^+$. Thus, by Lemma \ref{Uinv}, all fibers of
$p$ are reduced and irreducible. 

\begin{lemma}\label{vinberg}

\begin{enumerate}

\item The map $p:\cV\to \bA^{\Pi}$ is a family of reductive varieties
for $G$, the \emph{Vinberg family}.

\item The fiber of $p$ at any closed point of $\bA^{\Pi}$ with
non-zero coordinates is isomorphic to $G$. 

\item Let $X$ be an affine $G$-variety and let 
$$
p_X:X*\cV = (X\times \cV)//G \to \cV//G = \bA^{\Pi}
$$
be the family of $G$-varieties defined in Lemma \ref{construction}. 
Then the induced map $p_X//U:(X*\cV)//U\to \bA^{\Pi}$ is a trivial
family of $T$-varieties. As a consequence, if $X$ is irreducible,
then so are all fibers of $p_X$. Moreover, the fiber of $p_X$ at any
closed point with non-zero coordinates is isomorphic to $X$.

\end{enumerate}
\end{lemma}

\begin{proof}
(1) has just been proved.

(2) Note that the complement in $\bA^{\Pi}$ of the union of all
coordinate hyperplanes is a torus with character group the root
lattice. So this torus is the quotient of $T$ by the center $Z(G)$.
Moreover, the pull-back of $p$ to this torus is the natural map
$(T\times G)/Z(G)\to T/Z(G)$, where $Z(G)$ acts on $T\times G$ by
$z(t,g)=(zt,zg)$. Thus, its fibers are isomorphic to $G$.

(3) Since $X*\cV=(X\times\cV)//G$, one has
$$
k[X*\cV]^U = \bigoplus_{\lambda\in\Lambda^+,\,\mu\in\Lambda,\,
\lambda\leq_{\Pi}\mu} e^{\mu} k[X]^U_\lambda,
$$
a subalgebra of $k[T\times X]^U$. Thus, the $T$-invariant subspace 
$$
A=\bigoplus_{\lambda\in\Lambda^+} e^{\lambda} k[X]^U_\lambda
$$
is a subalgebra of $k[X*\cV]^U$, isomorphic to $k[X]^U$. Moreover,
$k[X*\cV]^U$ is isomorphic to $A\otimes_k R_0$ as a
$T$-algebra. This trivializes $p_X//U$, and proves the first
assertion. Together with Lemma \ref{Uinv}, this implies the second
assertion.
\end{proof}

Next we describe all fibers of $p_X$; for this, we recall a
construction that already occurred in the proof of Theorem 
\ref{Hilb}. Let $\gamma$ be a dominant one-parameter subgroup of 
$T$. For any non-negative integer $n$, let
$$
k[X]_{\le n}=
\bigoplus_{\lambda\in\Lambda^+,\,\langle\lambda,\gamma\rangle\le n} 
\Hom^G(V_\lambda,k[X])\otimes_k V_\lambda
$$
(where we identify each $\Hom^G(V_\lambda,k[X])\otimes_k V_\lambda$
with its image in $k[X]$ via the evaluation map).
Then the $k[X]_{\le n}$ form an ascending filtration of the
ring $k[X]$ by $G$-submodules. The associated graded ring 
$\gr_\gamma k[X]$ is isomorphic (as a $G$-module) to $k[X]$;
we may regard the product $V_\lambda\cdot V_\mu$ of any two simple
submodules of $\gr_\gamma k[X]$ as the subspace of their product in 
$k[X]$ spanned by those $V_\nu$ such that 
$\langle \gamma,\nu\rangle=\langle\gamma,\lambda+\mu\rangle$. 
Equivalently, $\nu\leq_K\lambda+\mu$, where $K$ is the set of simple
roots orthogonal to $\gamma$; in particular, $\gr_\gamma k[X]$
depends only on $K$. For example, $\gr_\gamma k[G]$ is the algebra
$k[G]_{(K)}$ considered in Subsection
\ref{subsec: Classification1}. 

We can now state the following result, whose proof is a direct checking.

\begin{lemma}\label{fibers}
Given a closed point $s=(s_\alpha)_{\alpha\in\Pi}\in \bA^{\Pi}$, let
$K=\{\alpha\in\Pi~\vert~s_\alpha=0\}$ and let $\gamma$ be a dominant
one-parameter subgroup of $T$, such that 
$K=\{\alpha\in\Pi~\vert~\langle\gamma,\alpha\rangle=0\}$.
Then the algebra of regular functions on the fiber of $p_X$ at $s$ is
isomorphic to $\gr_\gamma k[X]$.

In particular, the fiber of the Vinberg family at $s$ is the reductive
variety 
$X_{\Lambda^+_{\bR},K} = \Spec k[G]_{(K)} = 
\Spec k[\gxg/H_{K,\Lambda,\emptyset}]$.

\end{lemma}

If $s$ lies outside all coordinate hyperplanes of
$\bA^{\Pi}$, then $K=\emptyset$. Then we can take $\gamma=0$, whence
$(X * \cV)_s \simeq X$ (this also follows from Lemma
\ref{construction}). On the other hand, if $s$ is the origin of
$\bA^{\Pi}$, then $K=\Pi$, the dominant one-parameter subgroup
$\gamma$ is regular, and $\gr_\gamma k[X]$ is the 
degeneration of $k[X]$ constructed in \cite{Popov86} (see also
\cite{Grosshans97} Chapter 15).

\subsection{Local structure of families of reductive varieties}
\label{subsec: Local structure of families}

We show that the Vinberg family is a local model (in the sense of
Proposition \ref{curve}) for families of reductive varieties.

\begin{theorem}\label{locstruc}
Let $\pi:\cX\to S$ be a family of reductive varieties over an integral
scheme, with general fiber $X=X_\sigma$; let $K=K_\sigma$. Then every
$s\in S$ admits an open neighborhood $S'$ and a morphism 
$f:S'\to \bA^K$, such that the pull-back families 
$\cX\times_S S'$ and  $(X*\cV)\times_{\bA^{\Pi}}S'$ are isomorphic
(here we regard $\bA^K$ as a coordinate subspace of $\bA^{\Pi}$).
\end{theorem}

\begin{proof}
We may assume that $S$, and hence $\cX$, is affine, and also 
(by Lemma \ref{loctriv}) that the family 
$\pi//(U^-\times U):\cX//(U^-\times U)\to S$ is trivial. Let
$\cR=\Gamma(\cX,\cO_{\cX})$ and $\cR_0=\Gamma(S,\cO_S)$. As in the
proof of Corollary \ref{isotrivial}, write
$$
\cR = \bigoplus_{\lambda\in\Lambda^+\cap\sigma}
F_\lambda\otimes_k \End V_\lambda\subseteq
\bigoplus_{\lambda\in\Lambda^+\cap\sigma} 
k(\oS)\otimes_k\End V_\lambda,
$$ 
where every $F_\lambda$ is an invertible $\cR_0$-submodule of
$k(\oS)$. Then we may choose local generators
$f_\lambda$ of $F_\lambda$ at $s$, satisfying 
$f_\lambda f_\mu= f_{\lambda+\mu}$ for all $\lambda$, $\mu$ in
$\Lambda^+\cap\sigma$.

Since the group $\Lambda\cap\lin(\sigma)$ is generated by
$\Lambda^+\cap\sigma$, the assignment $\lambda\mapsto f_\lambda$
extends uniquely to a group homomorphism 
$\Lambda\cap\lin(\sigma)\to k(\oS)^*$. And since $K$ is
contained in $\Lambda\cap\lin(\sigma)$, we obtain in particular
elements $f_\alpha\in k(\oS)^*$ ($\alpha\in K$). 

Let $\alpha\in K$. By Lemma \ref{transvectants}, we may find 
$\lambda\in \Lambda^+\cap\sigma$ such that 
$\langle\lambda,\check\alpha\rangle\neq 0$; then the product 
$f_\lambda \End V_\lambda\cdot f_\lambda \End V_\lambda$
contains $f_{\lambda}^2 \End V_{2\lambda-\alpha}$. It follows that 
$$
f_{2\lambda} = f_\lambda^2 \in 
F_{2\lambda-\alpha} = f_{2\lambda-\alpha}\cR_0
= f_{2\lambda}f_{\alpha}^{-1}\cR_0,
$$ 
that is, $f_\alpha\in \cR_0$. Thus, we obtain a morphism
$f=(f_\alpha)_{\alpha\in K}:S\to\bA^K$.

Let $R_0$ be the subalgebra of $\cR_0$ generated by the $f_\alpha$,
$\alpha\in K$. We claim that the subspace of $\cR$:
$$
R=\bigoplus_{\lambda\in\Lambda^+\cap\sigma,\,\mu\in\Lambda,\,
\lambda\leq_K\mu} k f_\mu \otimes_k \End V_\lambda
=\bigoplus_{\lambda\in\Lambda^+\cap\sigma}
f_\lambda R_0\otimes_k \End V_\lambda
$$
is a subalgebra over $R_0$. Indeed, if 
$\lambda,\mu\in\Lambda^+\cap\sigma$ and $\nu\in\Lambda^+$ satisfies
$\nu\leq_K\lambda+\mu$, then 
$f_\lambda f_\mu = f_{\nu} f_{\lambda+\mu-\nu}\in f_\nu R_0$.

Clearly, the multiplication map
$$
R\otimes_{R_0} \cR_0\to \cR
$$ 
is an isomorphism of $\cR_0$-algebras. In other words, $\cX$ is the 
pull-back of the family $p_X: X*\cV\to\bA^{\Pi}$ under $f$.
\end{proof}

\providecommand{\bysame}{\leavevmode\hbox to3em{\hrulefill}\thinspace}



\end{document}